\documentclass{amsart}
\usepackage{amsmath,amsthm,amssymb,amscd}
\usepackage{wrapfig}

\renewcommand\appendix{\par
  \setcounter{section}{0}
  \setcounter{figure}{0}
  \setcounter{table}{0}
  \renewcommand\thesection{ \Alph{section}}
  \renewcommand\thefigure{\Alph{section}\arabic{figure}}
  \renewcommand\thetable{\Alph{section}\arabic{table}}
}

\newtheorem{thm}{Theorem}[section]
\newtheorem{prop}[thm]{Proposition}
\newtheorem{lem}[thm]{Lemma}
\newtheorem{cor}[thm]{Corollary}

\theoremstyle{definition}
\newtheorem{dfn}[thm]{Definition}
\theoremstyle{definition}
\newtheorem{eg}[thm]{Example}

\theoremstyle{remark}
\newtheorem{rem}[thm]{Remark}

\newcommand{\Z}{\mathbb{Z}}

\newcommand{\R}{\mathbb{R}}

\title{Fundamental groups of neighborhood complexes}

\begin{document}

\author{Takahiro Matsushita}
\email{mtst@math.kyoto-u.ac.jp}
\address{Department of Mathematics, Kyoto University, Oiwake-cho Sakyo-ku Kyoto 606-8502, Japan}

\begin{abstract}
The neighborhood complexes of graphs were introduced by Lov\'asz in his proof of the Kneser conjecture. He showed that a certain topological property of $N(G)$ gives a lower bound for the chromatic number of $G$.

In this paper, we study a combinatorial description of the fundamental groups of the neighborhood complexes. For a positive integer $r$, we introduce the $r$-fundamental group $\pi_1^r(G,v)$ of a based graph $(G,v)$ and the $r$-neighborhood complex $N_r(G)$ of $G$. The $1$-neighborhood complex is the neighborhood complex. We show that the even part $\pi_1^{2r}(G,v)_{ev}$, which is a subgroup of $\pi_1^{2r}(G,v)$ with index 1 or 2, is isomorphic to the fundamental group of $(N_r(G),v)$ if $v$ is not isolated. We can use the $r$-fundamental groups to show the non-existence of graph homomorphisms. For example, we show that $\pi_1^3(KG_{2k+1,k})$ is isomorphic to $\Z /2$, and this implies that there is no graph homomorphism from $KG_{2k+1,k}$ to the 5-cycle graph $C_5$. We discuss the covering maps associated to $r$-fundamental groups.
\end{abstract}

\maketitle

\section{Introduction}

An $n$-coloring of a graph $G$ is a map from the vertex set of $G$ to the $n$-point set $\{ 0,1,\cdots, n-1\}$ so that adjacent vertices have different values. The chromatic number of $G$ is the smallest integer $n$ such that $G$ has an $n$-coloring. The graph coloring problem, which is one of the most classical problems in graph theory, is to compute the chromatic number.

The neighborhood complex was introduced by Lov\'asz \cite{Lovasz} in the context of this subject. The neighborhood complex $N(G)$ of a graph $G$ is the simplicial complex whose vertices are non-isolated vertices and whose simplices are finite sets of vertices which have a common neighbor. Lov\'asz showed that if $N(G)$ is $n$-connected, then the chromatic number $\chi(G)$ of $G$ is greater than $n + 2$. Using this method, he determined the chromatic number of the Kneser graphs $KG_{n,k}$. For further development of the neighborhood complex and its related topics, we refer to \cite{Kozlov}.

In this paper, we study a combinatorial description of the fundamental group of the neighborhood complex. Throughout the paper, $r$ denotes a fixed positive integer. In Section 3, we introduce the $r$-fundamental group $\pi_1^r(G,v)$ of a based graph $(G,v)$ in a combinatorial way. If $G$ is simple, it turns out that $\pi_1^r(G,v)$ is isomorphic to the fundamental group of the 2-complex $X_r(G)$, which is obtained by attaching 2-cells to all $(2s)$-cycles of $G$ for $s \leq r$ (Theorem 3.1). The even part of $\pi_1^r(G,v) = \pi_1(X_r(G),v)$ is the subgroup of $\pi_1^r(G,v)$ consisting of elements represented by loops with even length, and is denoted by $\pi_1^r(G,v)_{ev}$.

In Section 4, we introduce the notion of the $r$-neighborhood complex $N_r(G)$. The 1-neighborhood complex is the neighborhood complex mentioned above. In terms of these notions, our principal result is formulated as follows.

\begin{thm}\label{Theorem 1.1}
Let $G$ be a graph, and let $v$ be a non-isolated vertex of $G$. Then there is a natural isomorphism
$$\begin{CD}
\pi_1(N_r(G),v) @>{\cong}>> \pi_1^{2r}(G,v)_{ev}.
\end{CD}$$
\end{thm}

We can apply the $r$-fundamental groups of graphs to the existence problem of the graph homomorphisms (see Section 2), which is a certain generalization of the graph coloring problem. In Section 6, we show that for every positive integer $k$, the 3-fundamental group of the Kneser graph $KG_{2k+1,k}$ is isomorphic to $\Z /2$ (Proposition \ref{Proposition 5.3}). This implies that there is no graph homomorphism from $KG_{2k+1, k}$ to the 5-cycle graph $C_5$ (Corollary \ref{Corollary 5.4}). Since the odd girth of $KG_{2k+1,k}$ is $2k+1$ (see Lemma \ref{Lemma 5.3.3}), this obstruction of the existence of graph homomorphisms is not obtained from the odd girths for $k \geq 2$. Moreover, this is not obtained from the topology of the neighborhood complex. In fact, the neighborhood complexes of $C_3$ and $C_5$ are homeomorphic to the circle $S^1$ and there is a graph homomorphism from $KG_{2k+1,k}$ to $C_3 \cong K_3$. The case of box complexes discussed in \cite{MZ} is the same.

In Section 6, we discuss the covering maps associated to the $r$-fundamental groups, called $r$-covering maps. As is the case of topological spaces (see Chapter 2 of \cite{Spanier}), it turns out that there is a 1-1 correspondence between subgroups of $r$-fundamental groups and connected $r$-coverings (Theorem \ref{Theorem 6.8}). Moreover, we have a purely graph theoretic description of covering spaces over neighborhood complexes. Namely, in case $G$ is non-bipartite, a connected covering space over $N_r(G)$ corresponds to a connected 2-covering over the Kronecker double covering $K_2 \times G$.


The rest of the paper is organized as follows. In Section 2, we review the definitions and facts related to graphs and simplicial complexes. In Section 3, we introduce the $r$-fundamental group $\pi_1^r(G,v)$ of a based graph $(G,v)$ and show that $\pi_1^r(G,v)$ is isomorphic to the fundamental group of a certain 2-dimensional complex $X_r(G)$. In Section 4, we define the $r$-neighborhood complex $N_r(G)$ and prove Theorem 1.1. In Section 5, we determine the $r$-fundamental groups of the Kneser graphs. In Section 6, we introduce the $r$-covering map and investigate its properties.

\section{Preliminaries}

In this section, we review definitions and facts related to graphs and simplicial complexes, in particular, the edge-path group of a based simplicial complex \cite{Spanier}. The edge-path group plays a central role in the proof of Theorem \ref{Theorem 1.1}.

A {\it graph} is a pair $G=(V(G),E(G))$ consisting of a set $V(G)$ together with a symmetric subset $E(G)$ of $V(G) \times V(G)$. We call an element of $V(G)$ a {\it vertex of $G$}. Hence our graphs are undirected, may have loops, but have no multiple edges. A graph is {\it simple} if it has no looped vertices. The {\it complete graph $K_n$ with $n$-vertices} is defined by $V(K_n) = \{ 0,1,\cdots, n-1\}$ and $E(K_n) = \{ (x,y) \; | \; x \neq y\}$. An $n$-coloring of $G$ is identified with a graph homomorphism from $G$ to $K_n$. Then the chromatic number of $G$ is formulated as the number
$$\chi(G) = \inf \{ n\geq 0 \; | \; \textrm{There is an $n$-coloring of $G$.}\}.$$
Here we consider the infimum of the empty set as $+\infty$.

An {\it (abstract) simplicial complex} is a pair $(V,\Delta)$ consisting of a set $V$ together with a family $\Delta$ of finite subsets of $V$, which satisfies the following conditions:
\begin{itemize}
\item $v \in V$ implies $\{ v\} \in \Delta$.
\item $\sigma \in \Delta$ and $\tau \subset \sigma$ imply $\tau \in \Delta$.
\end{itemize}
We call the set $V$ the {\it vertex set}, and an element of $\Delta$ a {\it simplex}. We frequently abbreviate $(V,\Delta)$ to $\Delta$. With this notation, we write $V(\Delta)$ to indicate the vertex set of $\Delta$. A {\it simplicial map from $\Delta_0$ to $\Delta_1$} is a map $f:V(\Delta_0) \rightarrow V(\Delta_1)$ such that $\sigma \in \Delta_0$ implies $f(\sigma) \in \Delta_1$.

Let $\Delta$ be a simplicial complex. The free $\R$-module generated by $V(\Delta)$ is denoted by $\R^{(V(\Delta))}$. We consider that the topology of $\R^{(V(\Delta))}$ is induced by the finite dimensional vector subspaces of $\R^{(V(\Delta))}$. For a vertex $v$ of $V(\Delta)$, the vector of $\R^{(V(\Delta))}$ associated to $v$ is denoted by $e_v$. For a simplex $\sigma$ of $\Delta$, let
$$\Delta_\sigma = \Big\{ \sum_{v \in \sigma} a_v e_v \; | \; a_v \geq 0, \; \sum_{v \in \sigma} a_v = 1 \Big\} \subset \R^{(V(\Delta))}.$$
The {\it geometric realization of $\Delta$} is defined to be the topological subspace
$$|\Delta| = \bigcup_{\sigma \in \Delta} \Delta_\sigma$$
of $\R^{(V(\Delta))}$. We frequently identify a vertex $v$ of $\Delta$ with the element $e_v$ of $|\Delta|$ associated to $v$.

We recall the definition of the edge-path group. An {\it edge-path} is a sequence $\gamma = (v_0,\cdots, v_n)$ of vertices of $\Delta$ such that $\{ v_{i-1} , v_i\}$ is a simplex of $\Delta$ for $i=1,\cdots, n$. We call $v_0$ the {\it origin of $\gamma$}, and $v_n$ the {\it end of $\gamma$}. Let $v$ and $w$ be vertices of $\Delta$. The set of edge-paths joining $v$ to $w$ is denoted by $\mathcal{E}(\Delta, v,w)$. Let $\gamma$ and $\gamma'$ be edge-paths of $\Delta$, and suppose that the end of $\gamma$ coincides with the origin of $\gamma'$. The {\it concatenation $\gamma' \cdot \gamma$ of edge-paths} is defined to be the sequence $\gamma$ followed by the sequence $\gamma'$, i.e. if $\gamma = (v_0,\cdots, v_m)$ and $\gamma' = (v'_0,\cdots, v'_n)$, then $\gamma' \cdot \gamma = (v_0,\cdots, v_n, v'_1,\cdots, v'_n)$.

Let $\sim$ be the equivalence relation on $\mathcal{E}(G,v,w)$ generated by the following two conditions:
\begin{itemize}
\item[(a)] If $v_{i-1}=v_i$, then
$$(v_0,\cdots, v_n) \sim (v_0,\cdots , v_{i-1}, v_{i+1}, \cdots, v_n).$$
\item[(b)] If $\{ v_{i-1},v_i,v_{i+1} \}$ is a simplex of $\Delta$, then
$$(v_0,\cdots, v_n) \sim (v_0,\cdots, v_{i-1},v_{i+1}, \cdots, v_n).$$
\end{itemize}
The edge-paths $\varphi$ and $\psi$ are {\it homotopic} if $\varphi \sim \psi$. Let $E(\Delta, v,w)$ be the quotient set $\mathcal{E}(\Delta,v,w) / \sim$.

A {\it based simplicial complex $(\Delta,v)$} is a simplicial complex $\Delta$ equipped with a distinct vertex $v$ of it. We write $E(\Delta,v)$ instead of $E(\Delta,v,v)$. Clearly the concatenation of the edge-paths induces a group operation of $E(\Delta,v)$, and we call this group the {\it edge-path group of $(\Delta,v)$}.

Let $\gamma = (v_0,\cdots, v_n)$ be an edge-path of $\Delta$. Define the path $\varphi_\gamma : [0,1] \rightarrow |\Delta|$ as follows:
\begin{itemize}
\item If $n = 0$, then we let $\varphi_\gamma : [0,1] \rightarrow |\Delta|$ be the constant path at $v_0$.
\item If $n > 0$, then
$$\varphi_\gamma \Big( \frac{i - 1 + t}{n}\Big) = (1-t)v_{i-1} + t v_i$$
for $i = 1,\cdots, n$ and $0 \leq t \leq 1$.
\end{itemize}
It is easy to show that the correspondence $\gamma \mapsto \varphi_\gamma$ induces a group homomorphism
$$\rho : E(\Delta,v) \longrightarrow \pi_1(|\Delta|,v), \; [\gamma] \mapsto [\varphi_\gamma].$$

\begin{prop}[Corollary 17 of Section 3.6 in \cite{Spanier}] \label{Proposition 2.1}
The group homomorphism $\rho$ is a natural isomorphism.
\end{prop}

\section{$r$-fundamental groups}

In this section we shall introduce the $r$-fundamental groups, and show that it is isomorphic to the fundamental group of a certain 2-dimensional cell complex $X_r(G)$ for a simple graph $G$. Recall that $r$ denotes a fixed positive integer.

For a non-negative integer $n$, let $L_n$ be the graph defined by $V(L_n) = \{ 0,1,\cdots, n\}$ and $E(L_n) = \{ (x,y)\; | \; |x-y| =1\}$. A graph homomorphism from $L_n$ to a graph $G$ is called a {\it path of $G$ with length $n$}. The length of a path $\varphi$ is denoted by $l(\varphi)$. We write $P(G,v,w)$ to indicate the set of paths joining $v$ to $w$. Consider the following conditions concerning a pair $\varphi$ and $\psi$ of paths joining $v$ to $w$:
\begin{itemize}
\item[(A)] $l(\psi) = l(\varphi) + 2$ and there is an integer $x \in \{ 0,1,\cdots, l(\varphi)\}$ such that $\psi(i) = \varphi(i)$ for $i \leq x$ and $\psi(i+2) = \varphi(i)$ for $i \geq x$ (see Figure 3.1).
\item[(B)$_r$] $l(\varphi) = l(\psi)$ and $\# \{ i \; | \; 0 \leq i \leq l(\varphi), \varphi(i) \neq \psi(i)\} < r$.
\end{itemize}
The equivalence relation generated by the conditions (A) and (B)$_r$ is denoted by $\simeq_r$. The paths $\varphi$ and $\psi$ are {\it $r$-homotopic} if $\varphi \simeq_r \psi$. Let $\pi_1^r(G,v,w)$ be the set of equivalence classes of $\simeq_r$ on $P(G,v,w)$. Let $[\varphi]_r$ denote the equivalence class of $\simeq_r$ represented by $\varphi$.

\begin{rem}\label{Remark 3.0}
Consider the following conditions concerning a pair $\varphi$ and $\psi$ of paths joining $v$ to $w$:
\begin{itemize}
\item[(B)$'_r$] $l(\varphi) = l(\psi)$ and there is $x \in \{ 0,1,\cdots, l(\varphi)\}$ such that $\varphi(i) = \psi(i)$ if $i \leq x$ or $i \geq x + r$.
\end{itemize}
Figure 3.2 describes the case $r = 3$. It is easy to see that the equivalence relation generated by (A) and (B)$'_r$ coincides with $\simeq_r$.
\end{rem}

\begin{center}
\begin{picture}(120,110)
\put(30,60){\circle*{3}} \put(70,50){\circle*{3}} \put(90,80){\circle*{3}} \put(100,20){\circle*{3}}

\put(0,40){\line(3,2){30}} \put(30,60){\line(4,-1){40}} \put(70,50){\line(2,3){20}} \put(70,50){\line(1,-1){30}} \put(100,20){\line(0,-1){20}}

\put(67,57){\vector(2,3){10}} \put(67,57){\line(-4,1){20}}
\put(78,51){\line(2,3){10}} \put(78,51){\vector(1,-1){15}}
\put(67,43){\vector(1,-1){17}} \put(67,43){\line(-4,1){23}}

\put(58,34){$\varphi$} \put(60,70){$\psi$} \put(95,80){$\psi(x+1)$}
\end{picture}

{\bf Figure 3.1}
\end{center}

\begin{center}
\begin{picture}(210,100)
\put(10,20){\circle*{3}} \put(50,10){\circle*{3}} \put(80,40){\circle*{3}} \put(120,50){\circle*{3}}
\put(140,80){\circle*{3}} \put(140,20){\circle*{3}} \put(170,20){\circle*{3}} \put(170,80){\circle*{3}} \put(190,50){\circle*{3}} \put(230,40){\circle*{3}}

\put(10,20){\line(-1,2){6}}
\put(10,20){\line(4,-1){40}} \put(50,10){\line(1,1){30}} \put(80,40){\line(4,1){40}} \put(120,50){\line(2,3){20}} \put(120,50){\line(2,-3){20}}
\put(140,80){\line(1,0){30}} \put(140,20){\line(1,0){30}} \put(170,20){\line(2,3){20}} \put(170,80){\line(2,-3){20}} \put(190,50){\line(4,-1){40}}
\put(230,40){\line(2,-3){10}}

\put(48,18){\vector(1,1){12}} \put(48,18){\line(-4,1){16}}

\put(115,56){\vector(2,3){8}} \put(115,56){\line(-4,-1){16}}
\put(117,42){\vector(2,-3){8}} \put(117,42){\line(-4,-1){15}}
\put(195,56){\line(-2,3){8}} \put(195,56){\vector(4,-1){16}}
\put(193,42){\vector(4,-1){15}} \put(193,42){\line(-2,-3){8}}

\put(128,49){$\varphi(x)$}
\put(118,75){$\psi$} \put(115,20){$\varphi$}
\end{picture}

{\bf Figure 3.2}
\end{center}

\vspace{2mm}
Let $\varphi : L_m \rightarrow G$ and $\psi : L_n \rightarrow G$ be paths and suppose $\varphi(m) = \psi(0)$. The {\it concatenation $\psi \cdot \varphi : L_{m+n} \rightarrow G$ of $\varphi$ followed by $\psi$} is defined by
$$\psi \cdot \varphi(i) = \begin{cases}
\varphi(i) & (i \leq m)\\
\psi(i - m) & (i \geq m).
\end{cases}$$
It is easy to see that the concatenation of paths induces a map
$$\pi_1^r(G , v,w) \times \pi_1^r(G,u,v) \longrightarrow \pi_1^r(G,u,w), ([\psi]_r, [\varphi]_r) \longmapsto [\psi \cdot \varphi]_r.$$

For a vertex $v$ of $G$, the graph homomorphism from $L_0$ to $G$ which takes $0$ to $v$ is denoted by $*_v$, and is called the {\it trivial loop at $v$}. For a path $\varphi : L_n \rightarrow G$, define $\overline{\varphi} : L_n \rightarrow G$ by $\overline{\varphi}(i) = \varphi(n - i)$ for $i = 0,1,\cdots, n$. For a path $\varphi$ joining $v$ to $w$, it is clear that $\varphi \cdot *_v = \varphi = *_w \cdot \varphi$ and $\overline{\varphi} \cdot \varphi \simeq_r *_v$.

A {\it based graph} is a pair $(G,v)$ consisting of a graph $G$ equipped with a distinct vertex $v$ of $G$. We call an element of $P(G,v,v)$ a {\it loop of $(G,v)$}. We write $\pi_1^r(G,v)$ instead of $\pi_1^r(G,v,v)$. The concatenation of loops induces a group operation of $\pi_1^r(G,v)$. In fact the identity element is $[*_v]$ and the inverse of $[\varphi]_r$ is $[\overline{\varphi}]_r$. We call this group the {\it $r$-fundamental group}. Note that a basepoint preserving graph homomorphism $f:(G,v) \rightarrow (H,w)$ induces a group homomorphism $f_* : \pi_1^r(G,v) \rightarrow \pi_1^r(H,w)$, $[\varphi]_r \mapsto [f \circ \varphi]_r$.

By the definition of the equivalence relation $\simeq_r$, we have a well-defined group homomorphism
$$l : \pi_1^r(G,v) \rightarrow \Z/2, ([\varphi]_r \mapsto (l(\varphi) \; {\rm mod}. 2)).$$
We call the kernel of this homomorphism the {\it even part of $\pi_1^r(G,v)$}, and write $\pi_1^r(G,v)_{ev}$ to indicate it. An element $\alpha$ of $\pi_1^r(G,v)$ is {\it even} if $\alpha$ is contained in $\pi_1^r(G,v)_{ev}$. An element of $\pi_1^r(G,v)$ is {\it odd} if it is not even. Note that the group homomorphism induced by a graph homomorphism preserves the parity. Clearly $\pi_1^r(G,v)_{ev}$ is a subgroup of $\pi_1^r(G,v)$ with index 1 or 2. The index of $\pi_1^r(G,v)_{ev}$ is $2$ if and only if the connected component of $G$ containing $v$ is not bipartite.

Let $v$ and $w$ be vertices of $G$, and let $\alpha$ be an $r$-homotopy class of paths joining $v$ to $w$. Then we have a group isomorphism
$$\pi_1^r(G,v) \rightarrow \pi_1^r(G,w), \; \beta \mapsto \alpha \cdot \beta \cdot \alpha^{-1}.$$
This isomorphism depends on the choice of $\alpha \in \pi_1^r(G,v,w)$. If a graph $G$ is connected and we only deal with the isomorphism class of the $r$-fundamental group, we often abbreviate $\pi_1^r(G,v)$ to $\pi^r_1(G)$.

Next we construct the 2-dimensional cell complex whose fundamental group is isomorphic to $\pi_1^r(G,v)$ in case $G$ is simple. Let $G$ be a simple graph. If we regard $G$ as a 1-dimensional simplicial complex in the usual way, we write $\Delta(G)$ to indicate the complex. We write $|G|$ instead of $|\Delta(G)|$.

Define the {\it $n$-cycle graph $C_n$} by $V(C_n) = \Z / n$ and $E(C_n) = \{ (x,x\pm 1) \; | \; x \in \Z/n\}$. An {\it $n$-cycle of $G$} is a graph homomorphism from $C_n$ to $G$. Here we do not require that the homomorphism be an embedding. We associate each $n$-cycle $C_n \rightarrow G$ to the continuous map $S^1 = |C_n| \rightarrow |G|$ in an obvious way. The 2-dimensional complex $X_r(G)$ is obtained from attaching 2-cells to $|G|$ by the maps $|C_{2s}| \rightarrow |G|$ associated to the $(2s)$-cycles of $G$ for $2 \leq s \leq r$.

\begin{thm} \label{Theorem 3.1}
Let $(G,v)$ be a based graph and suppose that $G$ is simple. Then there is a group isomorphism
$$\begin{CD}
\pi_1^r(G,v) @>{\cong}>> \pi_1(X_r(G),v).
\end{CD}$$
\end{thm}
\begin{proof}
Without loss of generality, we can assume that $G$ is connected.

First we define the group homomorphism $\overline{\Phi} : \pi_1^r(G,v) \rightarrow \pi_1(X_r(G),v)$. Let $\varphi : L_n \rightarrow G$ be a loop of $(G,v)$. We associate $\varphi$ to a path $\Phi(\varphi) : [0,1] \rightarrow X_r(G)$ of $X_r(G)$ as follows. If $n = 0$, then $\Phi(\varphi)$ is a constant path at $v \in X_r(G)$. If $n > 1$, let
$$\Phi( \varphi) \Big(\frac{k-1 +t}{n} \Big) =  (1-t)\varphi(k-1)+ t \varphi(k)$$
for $k= 1,\cdots, n$ and $0 \leq t \leq 1$. Here we consider the right side of the equality as an element of $|G| = |\Delta(G)|$. We write $\overline{\Phi}(\varphi)$ to indicate the homotopy class of the loop $\Phi(\varphi )$ of $(X_r(G),v)$.

Let $\varphi : L_m \rightarrow G$ and $\varphi' : L_n \rightarrow G$ be loops of $(G,v)$ and suppose $\varphi \simeq_r \varphi'$. We want to show that $\overline{\Phi}(\varphi) = \overline{\Phi}(\varphi')$. It suffices to show that $\Phi(\varphi) = \Phi(\varphi')$ in case the pair of loops $\varphi$ and $\varphi'$ satisfies one of the conditions (A) and (B)$'_r$ (see Remark 3.1). However, the verification of this part is straightforward, so we omit the details. Thus the correspondence $\varphi \mapsto \Phi (\varphi)$ induces a map $\overline{\Phi} : \pi_1^r(G,v) \rightarrow \pi_1(X_r(G),v)$. Clearly, $\overline{\Phi}$ is a group homomorphism.

Next we construct the inverse $\overline{\Psi} : \pi_1(X_r(G),v) \rightarrow \pi_1^r(G,v)$ of $\overline{\Phi}$. We start with the construction of a group homomorphism $\Psi : E(\Delta(G),v) \rightarrow \pi_1^r(G,v)$. Let $\alpha$ be an element of $E(\Delta(G),v)$. Let $\gamma = (v_0,\cdots, v_n)$ be an edge-path which represents $\alpha$ such that $v_i \neq v_{i-1}$ for $i=1,\cdots,k$. Let $\hat{\Psi}(\gamma) : L_n \rightarrow G$ be the loop which takes $i$ to $v_i$ for $i = 0,1,\cdots, n$. It is clear that if $\gamma \simeq \gamma'$, then $\hat{\Psi}(\gamma) \simeq_r \hat{\Psi}(\gamma')$. We let $\Psi(\alpha) = [\hat{\Psi}(\gamma)]_r$. Then $\Psi$ is a group homomorphism.

Let $N$ be the kernel of the composition of the sequence
$$\begin{CD}
E(\Delta(G),v) @>{\rho}>> \pi_1(|G|,v) @>{i_*}>> \pi_1(X_r(G),v).
\end{CD}$$
Here $i$ denotes the inclusion. By Proposition \ref{Proposition 2.1}, we have that $E(\Delta(G),v) / N \cong \pi_1(X_r(G),v)$. We now describe the generator of $N$. Let $\sigma$ be a $(2s)$-cycle $C_{2s} \rightarrow G$ with $s \leq r$. Let $\hat{\sigma}$ be the edge-path $(\sigma(0), \sigma(1), \cdots, \sigma(2s))$ of $\Delta(G)$. Let $\theta_\sigma$ be an edge-path joining $v$ to $\sigma(0)$. By Theorem 10 of Section 3.8 of \cite{Spanier} and Proposition \ref{Proposition 2.1}, $N$ is the normal subgroup of $E(\Delta(G),v)$ generated by the set $\{ \theta_\sigma \cdot \hat{\sigma} \cdot \overline{\theta}_\sigma \; | \; \textrm{$\sigma $ is a $(2s)$-cycle of $G$ with $s \leq r$.}\}$. Clearly, $\Psi : E(\Delta(G),v) \rightarrow \pi_1^r(G,v)$ takes an element of $N$ to $0$. Since $\pi_1(X_r(G),v) \cong E(\Delta(G),v) / N$, $\Psi$ induces a group homomorphism $\overline{\Psi} : \pi_1(X_r(G),v) \rightarrow \pi_1^r(G,v)$. It is clear that $\overline{\Psi}$ is the inverse of $\overline{\Phi}$. This completes the proof.
\end{proof}

Next we reduce the 2-cells attached to $|G|$. A $(2n)$-cycle $\sigma : C_{2n} \rightarrow G$ is {\it non-degenerate} if for $i \in \Z/{(2n)}$ and $k \in \Z/{(2n)}$, $\sigma (i) = \sigma(i + 2k)$ implies $k = 0$ modulo $n$. In case $G$ is bipartite, a $(2n)$-cycle of $G$ is non-degenerate if and only if it is embedding. Two $(2n)$-cycles $\sigma : C_{2n} \rightarrow G$ and $\tau : C_{2n} \rightarrow G$ are {\it equivalent} if there is an automorphism $f$ of $C_{2n}$ such that $\tau = \sigma \circ f$. The 2-dimensional complex $X'_r(G)$ is obtained by attaching 2-cells to $|G|$ by the maps associated to equivalence classes of non-degenerate $(2s)$-cycles for $2 \leq s \leq r$.

\begin{prop} \label{Proposition 3.3}
Let $(G,v)$ be a based graph and suppose that $G$ is simple. Then there is a group isomorphism
$$\begin{CD}
\pi_1^r(G,v) @>{\cong}>> \pi_1(X'_r(G),v).
\end{CD}$$
\end{prop}
\begin{proof}
For each equivalence class $\alpha$ of non-degenerate $(2s)$-cycles, we let $\sigma_\alpha$ be its representative. Regard $X'_r(G)$ as the subcomplex of $X_r(G)$, whose 2-cells are the 2-cells of $X_r(G)$ associated to $\sigma_\alpha$ for all $\alpha$. By Theorem 10 of Section 3.8 of \cite{Spanier}, it suffices to show that the attaching map of every 2-cell of $X_r(G)$ not contained in $X'_r(G)$ is null-homotopic in $X'_r(G)$.

Let $k$ be a positive integer with $k < r$ and suppose that every $l \leq k$ and every $(2l)$-cycle $\sigma : C_{2l} \rightarrow G$ of $G$, the map $S^1 \rightarrow |G|$ associated to $\sigma$ is null-homotopic in $|X'_r(G)|$. Let $\sigma' : C_{2k+2} \rightarrow G$ be a $(2k+2)$-cycle $\sigma'$. If $\sigma'$ is non-degenerate, then it follows from the definition of $X'_r(G)$ that the associated map $S^1 \rightarrow |G|$ is null-homotopic in $X'_r(G)$. Suppose that $\sigma'$ is degenerate. Then $\sigma'$ factors through $C_{2m} \vee C_{2m'}$ for $1 \leq m, m'  \leq k$. Here $C_{2m} \vee C_{2m'}$ denotes the wedge sum of $C_{2m}$ and $C_{2m'}$. It follows from the inductive hypothesis that the map associated to $\tau$ is null-homotopic.
\end{proof}

\begin{eg} \label{Example 3.4}
Let $n$ be a positive integer. The complex $X'_r(C_{2n})$ is a circle $S^1$ if $r < n$, and is a disk $D^2$ if $r \geq n$. Thus we have
$$\pi_1^r(C_{2n}) = \begin{cases}
\Z & (r < n)\\
1 & (r \geq n).
\end{cases}$$
On the other hand, $X'_r(C_{2n+1})$ is a circle if $r \leq 2n$ and is a 2-dimensional real projective plane if $r > 2n$. Thus we have
$$\pi_1^r(C_{2n+1}) = \begin{cases}
\Z & (r < n)\\
\Z / 2 & (r \geq n).
\end{cases}$$
\end{eg}

\section{Proof of Theorem 1.1}

In this section we give the definition of the $r$-neighborhood complex and prove Theorem 1.1.

The {\it neighborhood $N(v)$ of a vertex $v$ of $G$} is the set of vertices of $G$ adjacent to $v$. The {\it $r$-neighborhood $N_r(v)$} is defined inductively by $N_0(v) =  \{ v\}$ and
$$N_r(v) = \bigcup_{w \in N_{r-1}(v)} N(w).$$
The {\it $r$-neighborhood complex $N_r(G)$} is the simplicial complex defined as follows. The vertices of $N_r(G)$ are the non-isolated vertices of $G$, and a finite subset $\sigma$ of $V(G)$ is a simplex if and only if $\sigma$ is contained in the $r$-neighborhood of some vertex of $G$. The 1-neighborhood complex is the neighborhood complex introduced by Lov\'asz \cite{Lovasz}.

The rest of this section is devoted to the proof of Theorem \ref{Theorem 1.1}. It suffices to show that $E(N_r(G),v)$ and $\pi_1^{2r}(G,v)_{ev}$ are isomorphic (see Proposition \ref{Proposition 2.1}).

Now we construct a natural group homomorphism $\Phi : E(N_r(G),v) \rightarrow \pi_1^{2r}(G,v)_{ev}$. Let $\alpha$ be an element of $E(N_r(G),v)$ and let $\gamma = (x_0, \cdots, x_n)$ be an edge-path of $N_r(G)$ which represents $\alpha$. Then there is a loop $\varphi : L_{2rn} \rightarrow G$ of $(G,v)$ such that $\varphi(2ri) = x_i$ for $i= 0,1,\cdots ,n$. We set $\Phi(\alpha) = [\varphi]_{2r}$.

To see that $\Phi$ is well defined, we show that the $(2r)$-homotopy class $[\varphi]_{2r}$ does not depend on the choices of $\gamma$ and $\varphi$. First we show that for a fixed representative $\gamma$ of $\alpha$, the $(2r)$-homotopy class $[\varphi]_{2r}$ does not depend on the choice of $\varphi$. Consider another loop $\psi : L_{2rn} \rightarrow G$ of $(G,v)$ such that $\psi(2ri) = x_i$. Define the loop $\varphi_j$ $(j=0,1,\cdots,n)$ by the correspondence
$$\varphi_j(i) = \begin{cases}
\psi(i) & (i \leq 2 r j)\\
\varphi(i) & (i \geq 2rj).
\end{cases}$$
Then we have that $\varphi_0 = \varphi$, $\varphi_n = \psi$, and the pair of loops $\varphi_{j-1}$ and $\varphi_j$ satisfies the condition (B)$_{2r}$ in Section 3 for $j = 1,\cdots, n$. Hence we have $\varphi \simeq_{2r} \psi$. Therefore the $(2r)$-homotopy class $[\varphi]_{2r}$ does not depend on the choice of $\varphi$. So we shall write $\hat{\Phi}(\gamma)$ to indicate the $(2r)$-homotopy class $[\varphi]_{2r}$.

Next let $\gamma'$ be an edge-path of $N_r(G)$ which is homotopic to $\gamma$. We want to show that $\hat{\Phi}(\gamma) = \hat{\Phi} (\gamma')$. We can assume that the pair of edge-paths $\gamma$ and $\gamma'$ satisfies one of the conditions (a) and (b) in Section 2. The case of (a) is easily deduced, so we only give the proof of the case of (b). In this case $\gamma'$ is written as the edge-path $(x_0,\cdots, x_j, y, x_{j+1}, \cdots, x_n)$ such that the set $\{ x_j,y,x_{j+1}\}$ is a simplex of $N_r(G)$. By the definition of $N_r(G)$, there is a vertex $w$ of $G$ whose $r$-neighborhood contains $x_j, y$, and $x_{j+1}$. Let $\theta: L_{(2j+1)r} \rightarrow G$ be a path joining $x_0 = v$ to $w$ such that $\theta(2ri) = x_i$ for $i=0,1,\cdots, j$, and let $\theta' : L_{2(n-j)r - r} \rightarrow G$ be a path joining $w$ to $x_n = v$ such that $\theta'(2ri - r) = x_{j + i}$ for $i = 1,\cdots, n-j$. Finally, let $\tau: L_r \rightarrow G$ be a path joining $w$ to $y$. Consider the paths $\varphi = \theta' \cdot \theta$ and $\psi = \theta' \cdot \overline{\tau} \cdot \tau \cdot \theta$ (see Figure 4.1). It is clear that $\hat{\Phi}(\gamma) = [\varphi ]_{2r}$, $\hat{\Phi}(\gamma') = [\psi ]_{2r}$, and $\varphi$ and $\psi$ are $(2r)$-homotopic. Thus we have shown that the correspondence $\gamma \mapsto \hat{\Phi}(\gamma)$ induces the map $\Phi : E(N_r(G),v) \rightarrow \pi_1^{2r}(G,v)_{ev}$. It is clear that $\Phi$ is a group homomorphism, and natural with respect to basepoint preserving graph homomorphisms.

\begin{center}
\begin{picture}(320,160)(0,15)
\put(50,30){\circle*{3}} \put(80,70){\circle*{3}} \put(110,50){\circle*{3}} \put(150,70){\circle*{3}} \put(180,50){\circle*{3}}
\put(200,80){\circle*{3}} \put(240,70){\circle*{3}}

\put(50,30){\line(-2,1){20}} \put(50,30){\line(3,4){30}} \put(80,70){\line(3,-2){30}} \put(110,50){\line(2,1){40}} \put(150,70){\line(3,-2){30}} \put(180,50){\line(2,3){20}} \put(200,80){\line(4,-1){40}} \put(240,70){\line(1,-2){20}}

\put(170,100){\circle*{3}} \put(140,120){\circle*{3}} \put(150,150){\circle*{3}}

\put(150,70){\line(2,3){20}} \put(170,100){\line(-3,2){30}} \put(140,120){\line(1,3){10}}

\put(47,18){$x_j$} \put(145,59){$w$} \put(243,78){$x_{j+1}$} \put(157,153){$y$}

\put(145,76){\vector(2,3){12}} \put(145,76){\line(-2,-1){20}} \put(160,72){\vector(3,-2){18}} \put(160,72){\line(2,3){14}}

\put(150,55){\vector(3,-2){20}} \put(150,55){\line(-2,-1){24}} \put(165,28){$\varphi$} \put(142,90){$\psi$}
\end{picture}

{\bf Figure 4.1}
\end{center}

Next we construct the inverse $\Psi : \pi_1^{2r}(G,v)_{ev} \rightarrow E(N_r(G),v)$ of $\Phi$. By the hypothesis, there is a vertex $w$ adjacent to $v$. Let $\alpha$ be an element of $\pi_1^{2r}(G,v)_{ev}$. Let $\varphi : L_{2n} \rightarrow G$ be a loop of $(G,v)$ with even length. Let $m$ be an integer such that $2rm \geq n$. Define the extension $\varphi' : L_{2rm} \rightarrow G$ by
\begin{eqnarray}
\varphi'(i) = \begin{cases}
\varphi(i) & (i \leq 2n)\\
v & (i \textrm{ is even and } i\geq 2n)\\
w & (i \textrm{ is odd and }i \geq 2n).
\end{cases}
\end{eqnarray}
Define $\Psi(\varphi)$ to be the homotopy class of the edge-path
$$(\varphi'(0), \varphi'(2r), \cdots, \varphi'(2rm)).$$
Clearly the homotopy class $\Psi(\varphi)$ does not depend on the choice of the integer $m$.

We want to show that if the loops $\varphi$ and $\psi$ are $r$-homotopic, then the edge-paths $\Psi(\varphi)$ and $\Psi(\psi)$ are homotopic. We can assume that the pair of loops $\varphi$ and $\psi$ satisfies one of the conditions (A) and (B)$'_r$ in Section 3 (see Remark 3.1).

First we consider the case of (A). Namely, if $l(\varphi) = 2m$, then $l(\psi) = 2m+2$ and there is $x \in \{ 0,1,\cdots, 2m\}$ such that $\psi(i) = \varphi(i)$ for $i \leq x$ and $\psi (i + 2) = \varphi(i)$ for $i \geq x$. Let $m$ be an integer such that $2rm \geq 2n + 2$. Let $\varphi',\psi' : L_{2rm} \rightarrow G$ be the extensions of $\varphi, \psi$, respectively, described by the equation (1). Then the following claims hold:
\begin{itemize}
\item[(I)] $\{ \varphi' (2ri), \varphi'(2r(i+1)), \psi'(2r(i+1)) \} \subset N_r(\varphi'(2ri + r))$
\item[(II)] $\{ \varphi'(2ri), \psi'(2ri) , \psi'(2r(i+1))\} \subset N_r(\psi'(2ri + r))$
\end{itemize}
The claim (II) is obvious since $\varphi'(2ri)$ is equal to $\psi'(2ri)$ or $\psi'(2ri + 2)$. We shall show the claim (I). It is clear that $\varphi'(2ri)$ and $\varphi'(2r(i+1))$ are contained in $N_r(\varphi'(2ri + r))$. If $x \neq 2r(i+1) -1$, then $\psi'(2r(i+1))$ is equal to $\varphi'(2r(i+1))$ or $\varphi'(2r(i+1)-2)$, and is contained in $N_r(\varphi'(2ri+r))$. Suppose $x = 2r(i+1) -1$. Then $\psi'(x) = \varphi'(x) \in N_{r-1}(\varphi'(2ri+r))$, and hence $\psi'(2r(i+1)) = \psi'(x+1) \in N_r(\varphi'(2ri + r))$. This completes the proof of the claim (I).

By the above claims, the triangles appearing in Figure 4.2 form simplices of $N_r(G)$. Thus the edge-paths
$$(\varphi'(0), \varphi'(2r), \cdots, \varphi'(2rn))$$
and
$$(\psi'(0), \psi'(2r),\cdots, \psi'(2rn))$$
of $N_r(G)$ are homotopic.

\begin{center}
\begin{picture}(180,165)(0,0)
\put(40,30){\circle*{3}} \put(140,30){\circle*{3}} \put(40,130){\circle*{3}} \put(140,130){\circle*{3}}

\put(40,30){\line(1,0){100}} \put(40,30){\line(0,1){100}} \put(40,130){\line(1,0){100}} \put(140,30){\line(0,1){100}}
\put(40,130){\line(1,-1){100}}

\put(40,30){\line(-1,0){30}} \put(40,30){\line(-1,1){30}} \put(40,130){\line(-1,0){30}}
\put(140,30){\line(1,0){30}} \put(140,130){\line(1,0){30}} \put(140,130){\line(1,-1){30}}

\put(23,142){$\varphi'(2ri)$} \put(120,142){$\varphi'(2r(i+1))$}
\put(23,12){$\psi'(2ri)$} \put(120,12){$\psi'(2r(i+1))$}

\put(-10,128){$\cdots$} \put(180,128){$\cdots$}
\put(-10,28){$\cdots$} \put(180,28){$\cdots$}
\end{picture}

{\bf Figure 4.2}
\end{center}

Next we consider the case of (B)$'_r$ (see Remark 3.1). Then $l(\varphi) = l(\psi)$ and there is $x \in \{ 0,1,\cdots, l(\varphi)\}$ such that $\varphi(i) = \psi(i)$ if $i \leq x$ or $i \geq x+r$. Suppose that $x$ is even. Let $k$ be a non-negative integer such that $2r$ divides $2k + x$. Define $\varphi_1$ by
$$\varphi_1(i) = \begin{cases}
v  & (\textrm{$i$ is even and $i \leq 2k$.})\\
w & (\textrm{$i$ is odd and $i \leq 2k$}.)\\
\varphi(i - 2k) & ( i \geq 2k).
\end{cases}$$
Recall that $w$ is a fixed vertex adjacent to $v$. Define $\psi_1$ in a similar way. Then $\varphi_1(i) = \psi_1(i)$ if $i \leq 2k+x$ or $i \geq 2k+x + 2r$. Since $2r$ divides $2k+x$, we have that $\varphi_1(2ri) = \psi_1(2ri)$ for every integer $i$ with $0 \leq 2ri \leq l(\varphi) + 2k$. This implies that $\Psi(\varphi_1) = \Psi(\psi_1)$. Since we have already shown that $\Psi(\varphi)$ is invariant under the condition (A), we have $\Psi(\varphi) = \Psi(\varphi_1) = \Psi(\psi_1) = \Psi(\psi)$.

Next suppose that $x$ is odd. By the same way as the previous paragraph, we can suppose that $2r$ divides $x-1$. Let $\varphi',\psi' : L_{2rm} \rightarrow G$ be extensions of $\varphi$ and $\psi$ described in the equation (1), and suppose that $4r + x - 1 \leq 2rm$. Figure 4.3 illustrates $\varphi'$ and $\psi'$ in the case of $r = 2$. Then the following claims hold:
\begin{itemize}
\item[(i)] $\{ \varphi'(x-1), \varphi'(x+1), \varphi'(2r+x-1) \} \subset N_r(\varphi'(r+x-1))$
\item[(ii)] $\{ \varphi'(x-1), \varphi'(x+1), \psi'(2r + x- 1) \} \subset N_r(\psi'(r+x-1))$
\item[(iii)] $\{ \varphi'(x+1), \varphi'(2r+x-1), \psi'(2r+x-1)\} \subset N_r(\varphi'(r+x+1))$
\item[(iv)] $\{ \varphi'(2r+x -1), \psi'(2r+x-1), \varphi'(4r + x-1)\} \subset N_r(\varphi'(3r + x -1))$
\end{itemize}

Now we show the above claims. The claim (i) is obvious. The claim (ii) follows from $\varphi'(x) = \psi'(x) \in N_{r-1}(\psi'(r+x-1))$. The claim (iii) follows from $\psi'(2r+x) = \varphi'(2r+x) \in N_{r-1}(\varphi'(r+x+1))$. Finally, the claim (iv) follows from $\psi'(2r+x) = \varphi'(2r+x) \in N_{r-1}(\varphi'(3r+x-1))$.

Note that $\psi'(x-1) = \varphi'(x-1)$ and $\psi'(4r+x-1) = \varphi'(4r+x-1)$. It follows from the above claims and Figure 4.4 that the edge-paths $(\varphi'(x-1), \varphi'(2r+x-1), \varphi'(4r+x-1))$ and $(\psi'(x-1), \psi'(2r+x-1), \psi'(4r+x-1))$ are homotopic. Therefore we have that $(\varphi'(0),\varphi'(2r),\cdots, \varphi'(2rm))$ and $(\psi'(0),\psi'(2r),\cdots, \psi'(2rm))$ are homotopic, and hence $\Psi(\varphi) = \Psi(\psi)$. Thus the correspondence $\varphi \mapsto \Psi(\varphi)$ induces the map $\Psi : \pi_1^{2r}(G,v)_{ev} \rightarrow E(N_r(G),v)$. It is clear that $\Psi$ is the inverse of $\Phi$. This completes the proof of Theorem \ref{Theorem 1.1}.

\begin{center}
\begin{picture}(310,160)(0,10)
\put(20,80){\circle*{5}} \put(60,80){\circle*{3}}

\put(80,120){\circle*{5}} \put(120,140){\circle*{3}} \put(160,120){\circle*{5}}
\put(80,40){\circle*{3}} \put(120,20){\circle*{3}} \put(160,40){\circle*{5}}

\put(180,80){\circle*{3}} \put(220,60){\circle*{3}} \put(240,100){\circle*{3}} \put(280,80){\circle*{5}}

\put(20,80){\line(-2,1){23}} \put(280,80){\line(2,1){20}}

\put(20,80){\line(1,0){40}}
\put(60,80){\line(1,2){20}} \put(80,120){\line(2,1){40}} \put(120,140){\line(2,-1){40}} \put(160,120){\line(1,-2){20}}
\put(60,80){\line(1,-2){20}} \put(80,40){\line(2,-1){40}} \put(120,20){\line(2,1){40}} \put(160,40){\line(1,2){20}}

\put(180,80){\line(2,-1){40}} \put(220,60){\line(1,2){20}} \put(240,100){\line(2,-1){40}}


\put(55,88){\vector(1,2){8}} \put(55,88){\line(-1,0){16}} \put(46,99){$\varphi'$}
\put(55,72){\vector(1,-2){8}} \put(55,72){\line(-1,0){16}} \put(50,45){$\psi'$}

\put(186,86){\vector(2,-1){17}} \put(186,86){\line(-1,2){8}}
\put(183,70){\vector(2,-1){14}} \put(183,70){\line(-1,-2){8}}

\put(-10,63){$\varphi'(x-1)$} \put(33,128){$\varphi'(x + 1)$} \put(165,128){$\varphi'(2r+x-1)$}
\put(160,25){$\psi'(2r+x-1)$} \put(270,63){$\varphi'(4r+x-1)$}
\end{picture}

{\bf Figure 4.3}
\end{center}

\begin{center}
\begin{picture}(280,180)(0,-5)
\put(0,80){\circle*{3}} \put(80,80){\circle*{3}} \put(180,140){\circle*{3}} \put(180,20){\circle*{3}} \put(280,80){\circle*{3}}

\put(0,80){\line(1,0){80}} \put(0,80){\line(3,1){180}} \put(0,80){\line(3,-1){180}}
\put(80,80){\line(5,3){100}} \put(80,80){\line(5,-3){100}}
\put(180,20){\line(0,1){120}}
\put(280,80){\line(-5,3){100}} \put(280,80){\line(-5,-3){100}}

\put(-25,63){$\varphi'(x-1)$} \put(93,78){$\varphi'(x+1)$} \put(203,78){$\varphi'(4r+x-1)$}
\put(158,6){$\psi'(2r+x-1)$} \put(158,148){$\varphi'(2r+x-1)$}
\end{picture}

{\bf Figure 4.4}
\end{center}

\section{Kneser graphs}

Let $n$ and $k$ be positive integers, and suppose $n \geq 2k$. Define the Kneser graph $KG_{n,k}$ as follows. The vertices of $KG_{n,k}$ are the $k$-subsets of $\{ 0,1,\cdots, n-1\}$ and two of them are adjacent if and only if they are disjoint.

The Kneser conjecture asserts that the chromatic number of $KG_{n,k}$ is $n-2k+2$. In fact it is easy to show $\chi(KG_{n,k}) \leq n-2k+2$. On the other hand, Lov\'asz showed the following two theorems and solve the Kneser conjecture.

\begin{thm}[Lov\'asz \cite{Lovasz}] \label{Theorem 5.1}
If the neighborhood complex $N(G)$ of $G$ is $n$-connected, then $\chi(G) \geq n+3$.
\end{thm}

\begin{thm}[Lov\'asz \cite{Lovasz}] \label{Theorem 5.2}
The neighborhood complex of $KG_{n,k}$ is $(n-2k-1)$-connected.
\end{thm}

The purpose of this section is to determine the $r$-fundamental groups of the Kneser graphs for $r \geq 2$. If $n = 2k$, then the Kneser graph $KG_{2k,k}$ is the disjoint union of copies of $K_2$, so there is nothing we need to show.

Theorem \ref{Theorem 5.2} implies that $KG_{n,k}$ is connected if $n > 2k$. Suppose that $n > 2k+1$. Then Theorem \ref{Theorem 5.2} and Theorem \ref{Theorem 1.1} imply that the even part $\pi_1^2(KG_{n,k})_{ev}$ is trivial. Since the natural quotient map $\pi_1^2(G,v)_{ev} \rightarrow \pi_1^r(G,v)_{ev}$ is surjective, we have that $\pi_1^r(KG_{n,k})_{ev}$ is trivial for $r \geq 2$. Since $KG_{n,k}$ is not bipartite, the map
$$\pi_1^r(KG_{n,k}) \rightarrow \Z/2, \; [\gamma] \mapsto (l(\gamma) \; ({\rm mod.} \; 2))$$
is a surjective group homomorphism with trivial kernel $\pi_1^r(KG_{n,k})_{ev}$. Thus we have $\pi_1^r(KG_{n,k}) \cong \Z/2$.

Next we consider the case $KG_{2k+1,k}$. Note that $\pi_1^2(KG_{2k+1,k})$ is isomorphic to $\pi_1^1(KG_{2k+1,k})$. In fact, Lemma \ref{Lemma 5.3.2} mentioned below implies that every 4-cycle of $KG_{2k+1,k}$ is degenerate. Thus it follows from Proposition \ref{Proposition 3.3} that $\pi_1^2(KG_{2k+1,k}) \cong \pi_1^1(KG_{2k+1,k}) \cong \pi_1(|KG_{2k+1,k}|)$.

The number of the vertices of $KG_{2k+1,k}$ is
$$\binom{2k+1}{k}.$$
Since the degree of each vertex of $KG_{2k+1,k}$ is $k+1$, the number of edges is
$$\frac{k+1}{2} \binom{2k+1}{k}.$$
Thus the Euler characteristic $c(KG_{2k+1,k})$ of $KG_{2k+1,k}$ is
$$c(KG_{2k+1,k}) = \frac{1-k}{2} \binom{2k+1}{k}.$$
Let $F(m)$ denote the free group of rank $m$. For a connected graph $G$ with Euler characteristic $m$, it is clear that $\pi_1(G) \cong F(1-m)$. Thus we have
$$\pi_1^2(KG_{2k+1,k}) \cong F \bigg( 1 + \frac{k-1}{2} \binom{2k+1}{k} \bigg).$$
Next we compute the even part of $\pi_1^2(KG_{2k+1,k})$. Note that a subgroup of $F(m)$ with index 2 is isomorphic to $F(2m-1)$. To see this, let $X$ be the wedge of $m$ circles $S^1 \vee \cdots \vee S^1$. A subgroup of $F(m)$ with index 2 is isomorphic to some connected double covering $Y$ over $X$. Then the Euler characteristic $c(Y)$ of $Y$ is $2 c(X) = 2 - 2m$. Therefore we have
$$\pi_1(N(KG_{2k+1,k})) \cong \pi_1^2(KG_{2k+1,k})_{ev} \cong F\bigg( 1 + (k-1) \binom{2k+1}{k} \bigg).$$

Thus the only non-trivial part is the following proposition. It is clear that the proposition implies $\pi_1^r(KG_{2k+1,k}) \cong \Z/2$ for all $r \geq 3$.

\begin{prop} \label{Proposition 5.3}
The even part of $\pi_1^3(KG_{2k+1,k})$ is trivial. Hence $\pi_1^3(KG_{2k+1,k}) \cong \Z/2$.
\end{prop}

To prove this, we need to observe several properties of the Kneser graph $KG_{2k+1,k}$. For a non-negative integer $n$, we write $\langle n \rangle$ to indicate the set $\{ 0,1,\cdots, n-1\}$.

\begin{lem} \label{Lemma 5.3.1}
Let $\sigma$ and $\sigma'$ be $k$-subsets of $\langle 2k+1 \rangle$. There is a path of $KG_{2k+1,k}$ with length 2 joining $\sigma$ to $\sigma'$ in $KG_{2k+1,k}$ if and only if $\sigma$ and $\sigma'$ coincide except for one element.
\end{lem}
\begin{proof}
Let $\tau$ be a $k$-subset of $\langle 2k+1 \rangle$ adjacent to both $\sigma$ and $\sigma'$. Then $\sigma$ and $\sigma'$ are contained in the $(k+1)$-set $\langle 2k+1 \rangle \setminus \tau$. On the other hand, if $\sigma$ and $\sigma'$ coincide except for one element, then the cardinality of the union $\sigma \cup \sigma'$ is $k$ or $k+1$. Hence there is a $k$-subset $\tau$ of $\langle 2k+1 \rangle$ such that $\tau$ is disjoint from both $\sigma$ and $\sigma'$.
\end{proof}

\begin{lem} \label{Lemma 5.3.2}
Let $\sigma$ and $\sigma'$ be $k$-subsets of $\langle 2k+1 \rangle$ and suppose that $\sigma \neq \sigma'$ and there is a path with length 2 joining $\sigma$ to $\sigma'$. Then a vertex adjacent to both $\sigma$ and $\sigma'$ is unique.
\end{lem}
\begin{proof}
By Lemma \ref{Lemma 5.3.1}, the cardinality of the union $\sigma \cup \sigma'$ is $k+1$. Let $\tau$ be a $k$-subset of $\langle 2k+1 \rangle$ adjacent to both $\sigma$ and $\sigma'$. Then $\tau$ is contained in the $k$-subset $\langle 2k+1 \rangle \setminus (\sigma \cup \sigma')$. Since $\tau$ is a $k$-subset, we have $\tau = \langle 2k+1 \rangle \setminus (\sigma \cup \sigma')$.
\end{proof}

Recall that for a graph $G$, the {\it odd girth $g_0(G)$ of $G$} is the number
$$\inf \{ 2n+1  \; | \; \textrm{There is a graph homomorphism from $C_{2n+1}$ to $G$.}\}.$$
In other words, the odd girth is the minimal length of odd cycles embedded into $G$.

\begin{lem} \label{Lemma 5.3.3}
The odd girth of $KG_{2k+1,k}$ is $2k+1$.
\end{lem}
\begin{proof}
Let $\gamma : L_{2m} \rightarrow G$ be a graph homomorphism. Applying Lemma \ref{Lemma 5.3.1}, we can show that
$$\# (\gamma (0) \setminus \gamma(2m)) \leq m$$
by the induction on $m$. If $\gamma (0)$ and $\gamma(2m)$ are adjacent, we have that $\# (\gamma(0) \setminus \gamma(2m)) = \# \gamma(0) = k$, and hence $m \geq k$. Thus we have $g_0(KG_{2k+1,k}) \geq 2k+1$. On the other hand, let $\gamma_0 : \Z/(2k+1) \rightarrow V(KG_{2k+1,k})$ be the map defined by
$$\gamma_0(i) = \{ ki \; {\rm mod.} \; (2k+1), ki+1 \; {\rm mod.} \; (2k+1),\cdots, ki+k-1 \; {\rm mod.} \; (2k+1)\}.$$
Then $\gamma_0$ is a graph homomorphism from $C_{2k+1}$ to $KG_{2k+1,k}$. Thus we have $g_0(KG_{2k+1,k}) = 2k+1$.
\end{proof}

\begin{lem} \label{Lemma 5.3.4}
Let $\sigma$ and $\tau$ be $k$-subsets of $\langle 2k+1 \rangle$ and suppose that there is a path $\gamma$ with length 3 connecting $\sigma$ with $\tau$. Then the cardinality of the intersection $\sigma \cap \tau$ is at most 1.
\end{lem}
\begin{proof}
Let $\gamma : L_3 \rightarrow KG_{2k+1,k}$ be a path with length 3 joining $\sigma$ to $\tau$. By Lemma \ref{Lemma 5.3.1}, we have $\# (\sigma \setminus \gamma(2)) \leq 1$. Since $\sigma \cap \tau \subset \sigma \setminus \gamma(2)$, we have $\# (\sigma \cap \tau) \leq 1$.
\end{proof}

\begin{lem} \label{Lemma 5.3.5}
Let $\sigma$ and $\tau$ be vertices of $KG_{2k+1,k}$ and suppose that there is a path with length 3 joining $\sigma$ to $\tau$ and $\sigma \cap \tau \neq \emptyset$. In this case the cardinality of the union $\sigma \cup \tau$ is $2k-1$, and let $\{ a,b\} = \langle 2k+1 \rangle \setminus (\sigma \cup \tau)$. Then there are only two paths $\gamma$ and $\gamma'$ with length 3 joining $\sigma$ to $\tau$. Moreover, $a \not\in \gamma(1)$ implies $a \in \gamma'(1)$.
\end{lem}
\begin{proof}
It follows from Lemma \ref{Lemma 5.3.4} that $\sigma \cup \tau$ is a $(2k-1)$-subset of $\langle 2k+1 \rangle$. Let $\gamma : L_ 3 \rightarrow KG_{2k+1,k}$ be a path with length 3 joining $\sigma$ to $\tau$. Then $\gamma(1) \subset \langle 2k+1 \rangle \setminus \sigma$ and $\# (\gamma(1) \cap \tau) = k-1$ (see Lemma \ref{Lemma 5.3.1}). Since $\tau \setminus \sigma = k-1$ (see Lemma \ref{Lemma 5.3.4}), we have $\gamma(1) \cap \tau = \tau \setminus \sigma$ and hence $\tau \setminus \sigma \subset \gamma(1)$. Thus we have $\gamma(1) = (\tau \setminus \sigma) \cup \{ a\}$ or $(\tau \setminus \sigma) \cup \{ b\}$. If $\gamma(1)$ is determined, then $\gamma(2)$ is uniquely determined by Lemma \ref{Lemma 5.3.2}. Thus there exist only two paths joining $\sigma$ to $\tau$ with length 3. The last assertion clearly follows from this proof.
\end{proof}

Note that the two paths $\gamma$ and $\gamma'$ in Lemma \ref{Lemma 5.3.5} are $3$-homotopic.

Let $\sigma_0$ be the set $\langle k \rangle = \{ 0, 1, \cdots, k-1\}$ and consider $\sigma_0$ as the basepoint of $KG_{2k+1,k}$. Define the loop $\gamma_0 : L_{2k+1} \rightarrow KG_{2k+1,k}$ by
$$\gamma_0(i) = \{ ki + j \; {\rm mod.} \; (2k+1)\; | \; j \in \sigma_0\}.$$
It follows from Lemma \ref{Lemma 5.3.3} that the length of an odd loop is greater than or equal to $2k+1$.

\begin{lem} \label{Lemma 5.3.6}
Let $\alpha$ be an odd element of $\pi_1^3(KG_{2k+1,k},\sigma_0)$. Then there is a representative $\gamma$ of $\alpha$ whose length is minimal amongst the representatives of $\alpha$, such that $\gamma|_{L_{2k+1}} = \gamma_0$.
\end{lem}
\begin{proof}
Let $\gamma : L_{2m+1} \rightarrow KG_{2k+1,k}$ be a representative of $\alpha$ whose length is minimal amongst the representatives of $\alpha$. Suppose $\gamma(1) \neq \gamma_0(1)$. Then $\gamma(1)$ contains $2k$. Since $\gamma(2m+1) = \sigma_0$ does not contain $2k$, there is $a \in \{ 0,1,\cdots, m-1\}$ such that
$$\gamma(1),\; \gamma(3), \; \cdots, \; \gamma(2a+1)$$
contain $2k$ but $\gamma(2a+3)$ does not contain $2k$. Since the length of $\gamma$ is minimal, we have that $\gamma(2a) \cap \gamma(2a+3) \neq \emptyset$. In fact, if $\gamma(2a) \cap \gamma(2a+3) = \emptyset$, then the path $\hat{\gamma} : L_{2m-1} \rightarrow G$ defined by
$$\hat{\gamma}(i) = \begin{cases}
\gamma(i) & (i \leq 2a)\\
\gamma(i+2) & (i > 2a)
\end{cases}$$
is 3-homotopic to $\gamma$ (see Figure 5.1).

Thus there is a unique loop $\gamma' : L_{2m+1} \rightarrow G$ such that $\gamma'(i) = \gamma(i)$ if $i \leq 2a$ or $i \geq 2a+3$ but $\gamma \neq \gamma'$ (see Lemma 5.8). Then it is clear that $\gamma$ and $\gamma'$ are 3-homotopic. Since $\gamma(2a+1)$ contains $2k$, we have that $\gamma(2a)$ does not contain $2k$. Since neither $\gamma(2a)$ nor $\gamma(2a+3)$ contain $2k$, we have that $\gamma'(2a+1)$ does not contain $2k$ (see Lemma \ref{Lemma 5.3.5}). Thus
$$\gamma'(1), \; \gamma'(3), \; \cdots, \; \gamma'(2a-1)$$
contain $2k+1$ but $\gamma'(2a+1)$ does not contain $2k$. By the induction on $a$, there is a loop $\gamma'':L_{2m+1} \rightarrow KG_{2k+1,k}$ such that $\gamma \simeq_3 \gamma''$ and $\gamma''(1)$ does not contain $2k$. This implies that $\gamma''(1) = \gamma(1)$.

By the same argument, we can show that if $\gamma(i) = \gamma_0(i)$ $(i=1,\cdots, j)$ for some integer $j < 2k + 1$, then there is a loop $\gamma' : L_{2m+1} \rightarrow KG_{2k+1,k}$ such that $\gamma'(i) = \gamma_0(i)$ for $i = 0,1,\cdots, j+1$ and $\gamma' \simeq_3 \gamma$. This completes the proof.
\end{proof}

\begin{center}

\begin{picture}(120,80)(0,-10)
\put(40,20){\circle*{3}} \put(40,60){\circle*{3}} \put(80,20){\circle*{3}} \put(80,60){\circle*{3}}
\put(40,20){\line(1,0){40}} \put(40,60){\line(1,0){40}} \put(40,20){\line(0,1){40}} \put(80,20){\line(0,1){40}}
\put(40,20){\line(-2,-1){30}} \put(80,20){\line(2,-1){30}}

\put(32,25){\vector(0,1){20}} \put(32,25){\line(-2,-1){18}}
\put(88,25){\line(0,1){20}} \put(88,25){\vector(2,-1){18}}

\put(43,12){\line(-2,-1){15}} \put(43,12){\line(1,0){34}} \put(77,12){\vector(2,-1){15}}

\put(57,-3){$\hat{\gamma}$} \put(97,38){$\gamma$}
\end{picture}

{\bf Figure 5.1}

\end{center}

\begin{lem} \label{Lemma 5.3.7}
Let $(\sigma,\sigma')$ and $(\tau,\tau')$ be elements of $E(KG_{2k+1,k})$. Then there is an automorphism $\alpha$ of $KG_{2k+1,k}$ such that $\alpha(\sigma) = \tau$ and $\alpha(\sigma') = \tau'$.
\end{lem}
\begin{proof}
We write $S_{2k+1}$ to indicate the group of the automorphisms of the set $\langle 2k+1 \rangle = \{ 0,1,\cdots, 2k\}$. For an element $f \in S_{2k+1}$, define the automorphism $\alpha_f$ of $KG_{2k+1,k}$ by $\alpha_f(\sigma) = f(\sigma)$. Clearly, there is $f \in S_{2k+1}$ such that $\alpha_f(\sigma) = \tau$ and $\alpha_f(\sigma') = \tau'$ since $\sigma \cap \sigma' = \emptyset$ and $\tau \cap \tau' = \emptyset$.
\end{proof}

\noindent
{\it Proof of Proposition 5.3.} Let $\alpha$ be an odd element of $\pi_1^3(KG_{2k+1,k})$, let $\gamma$ be a representative of $\alpha$ whose length is minimal amongst the representatives of $\alpha$. We want to show that $\gamma$ is $3$-homotopic to $\gamma_0$. Since the odd girth of $KG_{2k+1,k}$ is $2k+1$, we have $l(\gamma) \geq 2k+1$. If $l(\gamma) = 2k+1$, then Lemma \ref{Lemma 5.3.6} implies $\gamma \simeq_3 \gamma_0$. Suppose that $l(\gamma)$ is greater than $2k+1$. By Lemma \ref{Lemma 5.3.6} we can assume that $\gamma|_{L_{2k+1}} = \gamma_0$, and set $\gamma = \gamma_1 \cdot \gamma_0$.

It follows from Lemma \ref{Lemma 5.3.7} that there is an automorphism $\alpha$ of $KG_{2k+1,k}$ such that $\alpha(\sigma_0) = \sigma_0$ and $\alpha(\gamma_0(2k)) = \gamma_1(1)$. Let $\gamma'_0 = \alpha \circ \gamma_0$. Then we have $\gamma_1 \cdot \gamma_0 \simeq_3 \gamma_1 \cdot \gamma'_0$ since $\gamma_0 \simeq_3 \gamma'_0$ (see the previous paragraph).

Note that $\gamma_1 \cdot \gamma'_0 (2k) = \gamma'_0(2k) = \alpha(\gamma_0(2k)) = \gamma_1(1) = \gamma_1 \cdot \gamma'_0 (2k+2)$. Hence $\gamma$ is 3-homotopic to a loop whose length is smaller than $l(\gamma)$ (see the condition (A) in Section 3). This contradicts the assumption of $\gamma$. Thus we have $l(\gamma) = 2k+1$ and $\gamma \simeq_3 \gamma_0$.

Therefore we conclude that $\pi_1^3(KG_{2k+1,k})$ has only one odd element. Since the index of the even part is 2, this implies that the even part of $\pi_1^3(KG_{2k+1,k})$ is trivial.\qed

\vspace{2mm}
We conclude this section with the following corollary.

\begin{cor}\label{Corollary 5.4}
For every positive integer $k$, there is no graph homomorphism from $KG_{2k+1,k}$ to $C_5$.
\end{cor}
\begin{proof}
Suppose that there is a graph homomorphism $f: KG_{2k+1,k} \rightarrow C_5$. Since $\pi_1^3(KG_{2k+1,k}) \cong \Z/2$ and $\pi_1^3(C_5) \cong \Z$ (see Example \ref{Example 3.4}), the group homomorphism $f_* : \pi_1^3(KG_{2k+1,k}) \rightarrow \pi_1^3(C_5)$ induced by $f$ is trivial. Let $\alpha$ be the generator of $\pi_1^3(KG_{2k+1,k})$. Since $\alpha$ is odd, we have that $f_* (\alpha)$ is odd and hence non-trivial. This is a contradiction.
\end{proof}

Note that if there is a graph homomorphism $f:G \rightarrow H$, then we have $g_0(G) \geq g_0(H)$. Since the odd girth of $KG_{2k+1,k}$ is $2k+1$ (Lemma \ref{Lemma 5.3.3}), the obstruction of the existence of graph homomorphisms from $KG_{2k+1,k}$ to $C_5$ is not obtained from the odd girths if $k \geq 2$. Moreover, we should note that this obstruction is not obtained from the topology of neighborhood complexes. In fact, both $N(C_5)$ and $N(C_3)$ are homeomorphic to $S^1$ but there is a graph homomorphism from $KG_{2k+1,k}$ to $C_3 \cong K_3$ (recall $\chi(KG_{n,k}) = n-2k+2$). Similarly, this obstruction is not obtained from the equivariant topology of box complexes discussed in \cite{MZ}.

\section{Covering maps}

It is well known that there is a close relationship between covering spaces and fundamental groups of topological spaces. In this section we introduce the covering notion associated to $r$-fundamental groups, called $r$-covering maps, and investigate their properties. As is the case of topogical spaces, there is a correspondence between a subgroup of $\pi_1^r(G,v)$ and an connected $r$-covering over $(G,v)$ (Theorem \ref{Theorem 6.8}).

It turns out that the even part of $\pi_1^r(G,v)$ corresponds to the Kronecker double covering $K_2 \times G$ (Proposition \ref{Proposition Kronecker}). Since the fundamental group of the $r$-neighborhood complex is isomorphic to the even part of $(2r)$-fundamental group, this implies that there is a connected covering space over $N_r(G)$ and a $(2r)$-covering over $K_2 \times G$. We also show that $(2r)$-covering induces a covering between $r$-neighborhood complexes (Proposition \ref{Proposition 6.20}).

Now we start with the definition of the $r$-covering maps.

\begin{dfn}
A graph homomorphism $p:G \rightarrow H$ is an {\it $r$-covering map} if the map
$$p|_{N_i(v)} : N_i(v) \rightarrow N_i(p(v)) $$
is a bijection for every $1 \leq i \leq r$ and every vertex $v$ of $G$.
\end{dfn}

\begin{lem} \label{Lemma 6.1}
Let $p:G \rightarrow H$ be a graph homomorphism. Then $p$ is an $r$-covering map if and only if for every vertex $v$ of $G$, the map
$$p|_{N(v)} : N(v) \rightarrow N(p(v))$$
is surjective and the map
$$p|_{N_r(v)} : N_r(v) \rightarrow N_r(p(v))$$
is injective.
\end{lem}
\begin{proof}
Since the ``only if" part is obvious, we only show the ``if" part. Suppose that $p|_{N(v)}$ is surjective and $p|_{N_r(v)}$ is injective for every $v \in V(G)$. It is straightforward to show that $p|_{N_i(v)}$ is surjective and $i \geq 0$ by the induction on $i$, and we omit the details. We show the injectivity of $p|_{N_i(v)}$ for every $v \in V(G)$ and $i = 1,\cdots ,r$. In case $v$ is isolated, the vertex $p(v)$ of $H$ is also isolated since $p|_{N(v)} : N(v) \rightarrow N(p(v))$ is surjective. Hence we have that $N_i(v) = \emptyset$ and $N_i(p(v)) = \emptyset$, and hence $p|_{N_i(v)} : N_i(v) \rightarrow N_i(p(v))$ is injective. Suppose that $v$ is not isolated. Then there is an element $w$ of $N_{r-i}(v)$. Since $N_i(v) \subset N_r(w)$, the injectivity of $p|_{N_i(v)}$ follows from the injectivity of $p|_{N_r(w)}$.
\end{proof}

\begin{eg} \label{Example 6.1.1}
Let $n$ be a positive integer, and let $m$ be an integer greater than 1. Then the graph homomorphism
$$p : C_{2nm} \rightarrow C_{2n}, \; (x \; {\rm mod.} \; 2nm) \mapsto (x \; {\rm mod.} \; 2n)$$
is an $(n-1)$-covering but not an $n$-covering. On the other hand, the graph homomorphism $p: C_{2(2n+1)} \rightarrow C_{2n+1}$ is an $r$-covering map for every positive integer $r$. If $m > 2$, then the graph homomorphism $C_{m (2n+1)} \rightarrow C_{2n+1}$ is a $(2n)$-covering but not a $(2n+1)$-covering. Compare these examples with the $r$-fundamental groups of cycles (Example \ref{Example 3.4}).
\end{eg}

\begin{lem} \label{Lemma 6.2}
Let $p:G\rightarrow H$ and $q: H \rightarrow K$ be graph homomorphisms. Suppose that $p:V(G) \rightarrow V(H)$ is surjective as a set map. If two of $p$, $q$, and $q \circ p$ are $r$-covering maps, then so is the third.
\end{lem}
\begin{proof}
For a vertex $v$ of $G$ and $i=1,\cdots, r$, consider the diagram
$$\begin{CD}
N_i(v) @>p>> N_i(p(v))\\
@V{q \circ p}VV @VVqV\\
N_i(q\circ p(v)) @= N_i(q \circ p(v)).
\end{CD}$$

\vspace{1mm} \noindent If two of the three arrows are bijective, then so is the third.
\end{proof}

\begin{prop} \label{Proposition 6.3}
Let $p:G \rightarrow H$ be an $r$-covering map and let $v$ be a vertex of $G$. Then the following hold.

\begin{itemize}
\item[(1)] Let $\varphi : L_m \rightarrow H$ be a path of $H$ starting from $p(v)$. Then there is a unique path $\tilde{\varphi}$ of $G$ such that $p \circ \tilde{\varphi} = \varphi$ and $\tilde{\varphi}(0) = v$. We call this $\tilde{\varphi}$ the lift of $\varphi$ starting from $v$.
\item[(2)] Let $\varphi : L_m \rightarrow H$ and $\psi : L_n \rightarrow H$ be paths such that $\varphi(0) = \psi(0)$ and $\varphi(m) = \psi (n)$. Let $\tilde{\varphi}$ and $\tilde{\psi}$ be the lifts of $\varphi$ and $\psi$ starting from $v$ respectively. If $\varphi \simeq_r \psi$ then $\tilde{\varphi} (m) = \tilde{\psi}(n)$ and $\tilde{\varphi} \simeq_r \tilde{\psi}$.
\end{itemize}
\end{prop}
\begin{proof}
The proof of (1) is straightforward and is omitted.

Now we show (2). We can assume that the pair of paths $\varphi: L_m \rightarrow H$ and $\psi : L_n \rightarrow H$ satisfies one of the conditions (A) and (B)$'_r$ in Section 3. Suppose that the condition (A) holds. Then $n = m+2$ and there is $x \in \{ 0,1,\cdots, m\}$ such that $\varphi(i) = \psi(i)$ for every $i \leq x$, and $\varphi(i) = \psi(i+2)$ for every $i \geq x$. Note that $\tilde{\psi}(x), \tilde{\psi} (x+2) \in N(\tilde{\psi}(x+1))$ and
$$p(\tilde{\psi}(x)) = \psi(x) = \varphi(x) = \psi(x+2) = p(\tilde{\psi}(x+2)).$$
Since $p$ is an $r$-covering map, we have $\tilde{\psi}(x) = \tilde{\psi}(x+2)$. The uniqueness of the lift (see (1)) implies $\tilde{\psi}(i) = \tilde{\varphi}(i)$ for every $i \leq x$. Hence we have $\tilde{\psi}(x+2) = \tilde{\psi}(x) = \tilde{\varphi}(x)$. The uniqueness of the lift again implies $\tilde{\psi}(i+2) = \tilde{\varphi}(i)$ for every $i \geq x$. Hence the pair $\tilde{\varphi}$ and $\tilde{\psi}$ of paths of $G$ satisfies the condition (A).

Next suppose that the pair $\varphi$ and $\psi$ satisfies the condition (B)$'_r$, namely, $m=n$ and there is $x \in \{ 0,1,\cdots, n\}$ such that $\varphi(i) = \psi (i)$ if either $i \leq x$ or $i \geq x+r$ holds. By the uniqueness of the lift, we have that $\tilde{\varphi}(i) = \tilde{\psi}(i)$ for $i \leq x$. Suppose $x+r \leq n$. Since $\tilde{\varphi}(x) = \tilde{\psi}(x)$, we have $\tilde{\varphi}(x+r), \tilde{\psi}(x+r) \in N_r(\tilde{\varphi}(x))$. Moreover, we have
$$p ( \tilde{\psi}(x+r)) = \psi(x+r) = \varphi(x+r) = p ( \tilde{\varphi}(x+r)).$$
Since $p$ is an $r$-covering map, we have $\tilde{\varphi}(x+r) = \tilde{\psi}(x+r)$. Again the lift of the uniqueness implies that $\tilde{\varphi}(i) = \tilde{\psi}(i)$ for every $i \geq x+r$. Thus the pair $\tilde{\varphi}$ and $\tilde{\psi}$ of paths of $G$ satisfies the condition (B)$'_r$. Since the proof in the case $x + r > n$ is similar, we omit it.
\end{proof}

\begin{cor} \label{Corollary 6.4}
Let $p:G \rightarrow H$ be an $r$-covering map, and let $v$ be a vertex of $G$. Then the group homomorphism $p_* : \pi_1^r(G,v) \rightarrow \pi_1^r(H,p(v))$ induced by $p$ is injective.
\end{cor}
\begin{proof}
Let $\varphi$ be a loop of $(G,v)$ and suppose that $p \circ \varphi$ is $r$-homotopic to $*_{p(v)}$. It follows from (2) of Proposition \ref{Proposition 6.3} that $\varphi$ is $r$-homotopic to $*_v$.
\end{proof}

\begin{cor} \label{Corollary 6.4.1}
Let $p:(G,v) \rightarrow (H,w)$ be a basepoint preserving $r$-covering map, and let $\varphi : L_n \rightarrow H$ be a loop of $(H,w)$. Then $[\varphi]_r$ belongs to $p_* (\pi_1^r(G,v))$ if and only if the lift $\tilde{\varphi}$ of $\varphi$ starting from $v$ is a loop of $(G,v)$.
\end{cor}
\begin{proof}
Suppose that $[\varphi]_r$ belongs to $p_* (\pi_1^r(G,v))$. Let $\psi$ be a loop of $(G,v)$ such that $p \circ \psi \simeq_r \varphi$. By (2) of Proposition \ref{Proposition 6.3}, we have that $\tilde{\varphi}$ is a loop.
\end{proof}

\begin{lem} \label{Lemma 6.5}
Let $p:(G,v) \rightarrow (H,w)$ be a basepoint preserving $r$-covering map, let $(T,x)$ be a connected based graph, and let $f:(T,x) \rightarrow (H,w)$ be a basepoint preserving graph homomorphism. Then a basepoint preserving graph homomorphism $\tilde{f}:(T,x) \rightarrow (G,v)$ such that $p \circ \tilde{f} = f$ exists if and only if $f_* (\pi_1^r(T,x)) \subset p_*(\pi_1^r(G,v))$. Moreover, such a graph homomorphism $\tilde{f}$ is uniquely determined.
\end{lem}
\begin{proof}
If there is a lift $\tilde{f}$ of $f$, then we have $f_* (\pi_1^r(T,x)) = p_* \circ \tilde{f}_* (\pi_1(T,x)) \subset p_* (\pi_1^r(G,v))$. On the other hand, suppose that $f_*(\pi_1^r(T,x))$ is contained in $p_* (\pi_1^r(G,v))$. We construct a graph homomorphim $\tilde{f}$ as follows. Let $y \in V(T)$ and let $\varphi : L_n \rightarrow T$ be a path joining $x$ to $y$. Let $\tilde{\varphi} : L_n \rightarrow G$ be the lift of $f \circ \varphi$ starting from $v$.

We show that the terminal point of $\tilde{\varphi}$ does not depend on the choice of the path $\varphi$. Let $\psi$ be another path joining $x$ to $y$. Let $\gamma$ be the lift of $f \circ (\overline{\psi} \cdot \varphi)$ of $G$ whose initial point is $v$. By the hypothesis $f_* (\pi_1^r(T,x)) \subset p_* (\pi_1^r(G,v)) $ and Corollary \ref{Corollary 6.4.1} we have that $\gamma$ is a loop of $(G,v)$. Hence the terminal point of $\tilde{\psi}$ coincides with the one of $\tilde{\psi} \cdot \gamma $. Note that $\tilde{\psi} \cdot \gamma$ is the lift of $p \circ (\psi \cdot \overline{\psi} \cdot \varphi )$. Since $\psi \cdot \overline{\psi} \cdot \varphi \simeq_r \varphi$, we have that the terminal points of $\tilde{\psi} \cdot \gamma$ and $\tilde{\varphi}$ coincide (see (2) of Proposition \ref{Proposition 6.3}).

Thus we let $\tilde{f}(y)$ be the terminal point of $\tilde{\varphi}$. It is straightforward to show that $\tilde{f}$ is a basepoint preserving graph homomorphism from $(T,x)$ to $(G,v)$ and $p \circ \tilde{f} = f$, and we omit the details.

Finally, we show that the uniqueness of the lift $\tilde{f}$. Let $\tilde{f}_0 : (T,x) \rightarrow (G,v)$ be a basepoint preserving graph homomorphism such that $p \circ \tilde{f}_0 = f$. Consider the set $A = \{ y \in V(T) \; | \; \tilde{f}_0(y) = \tilde{f}(y)\}$. Since $\tilde{f}_0(x) = v = \tilde{f}(x)$, we have that $x$ belongs to $A$. Suppose $y \in A$ and let $z \in N(y)$. Note that $\tilde{f}(z)$ and $\tilde{f}_0(z)$ are contained in $N(\tilde{f}(y))$ and $p(\tilde{f}_0(z)) = f(z) = p(\tilde{f}(z))$. Since $p$ is an $r$-covering map, we have $\tilde{f}_0(z) = \tilde{f}(z)$, and hence $z \in A$. Since $T$ is connected, the set $A$ coincides with the vertex set of $T$. This implies $\tilde{f}_0 = \tilde{f}$.
\end{proof}

The {\it universal $r$-covering over a based graph $(G,v)$} is a basepoint preserving $r$-covering map $p:(\tilde{G}, \tilde{v}) \rightarrow (G,v)$ such that $\tilde{G}$ is connected and the $r$-fundamental group of $(\tilde{G}, \tilde{v})$ is trivial. It follows from Lemma \ref{Lemma 6.5} that the universal $r$-covering is unique up to isomorphisms.

\begin{prop} \label{Proposition 6.6}
For every based graph $(G,v)$, the universal $r$-covering over $(G,v)$ exists.
\end{prop}
\begin{proof}
Recall that $\pi_1^r(G,v,w)$ denotes the set of $r$-homotopy classes of paths joining $v$ to $w$ (see Section 3). Set
$$V(\tilde{G})=\coprod_{w \in V(G)} \pi_1^r(G,v,w),$$
$$E(\tilde{G}) = \{ (\alpha,\beta) \; | \; \textrm{There is $\varphi \in \beta $ such that the path $\varphi |_{L_{l(\varphi)-1}}$ belongs to $\alpha $.}\},$$
and let $p:V(\tilde{G}) \rightarrow V(G)$ be the map which takes an element of $\pi_1^r(G,v,w)$ to $w$. Let $\tilde{v}$ be the $r$-homotopy class of the trivial loop $*_v$ of $(G,v)$.

First we show that $E(\tilde{G})$ is a symmetric subset of $V(\tilde{G}) \times V(\tilde{G})$. Consider the following conditions (1), (2), and (3) concerning a pair $(\alpha,\beta)$ of elements of $V(\tilde{G})$.
\begin{itemize}
\item[(1)] The pair $(\alpha,\beta)$ belongs to $E(\tilde{G})$.
\item[(2)] There is a representative $\varphi$ of $\alpha$ such that the map $\varphi ' : V(L_{l(\varphi) + 1}) \rightarrow V(G)$ defined by $\varphi'|_{V(L_{l(\varphi)})} = \varphi$ and $\varphi' (l(\varphi) + 1) = p(\beta)$ is a graph homomorphism, and $\varphi '$ belongs to $\beta$.
\item[(3)] For each representative $\varphi$ of $\alpha$, the map $\varphi ' :V(L_{l(\varphi) + 1}) \rightarrow V(G)$ defined by $\varphi' |_{V(L_{l(\varphi)})} = \varphi$ and $\varphi' (l(\varphi) + 1) = p(\beta)$ is a graph homomorphism. Moreover, the path $\varphi'$ belongs to $\beta$.
\end{itemize}
We show that these conditions are equivalent. It is clear that the conditions (1) and (2) are equivalent and the condition (3) implies the condition (2). Suppose that the condition (2) holds. Let $u : L_1 \rightarrow G$ be a path which takes $0$ to $p(\alpha)$ and $1$ to $p(\beta)$. Then there is a representative $\varphi$ of $\alpha$ such that the extension $\varphi' = u \cdot \varphi : L_{l(\varphi) + 1} \rightarrow G$ is a representative of $\beta$. Let $\psi$ be another representative of $\alpha$ and consider $\psi' = u \cdot \psi$. Since $\varphi \simeq_r \psi$, we have that $\psi'  = u \cdot \psi \simeq_r u \cdot \varphi = \varphi' \in \beta$. Thus the condition (3) holds.

Let $(\alpha,\beta) \in E(\tilde{G})$. Let $\varphi$ be an element of $\alpha$. Define the path $\varphi'' : L_{l(\varphi) + 2} \rightarrow G$ by $\varphi'' |_{L_{l(\varphi)}} = \varphi$, $\varphi'' (l(\varphi) +1) = p(\beta)$, and $\varphi''(l(\varphi) +2) = p(\alpha)$. Then we have $(\beta,\alpha) = ([\varphi']_r , [\varphi]_r) = ([\varphi']_r, [\varphi'']_r) \in E(\tilde{G})$. Hence $E(\tilde{G})$ is symmetric.

Next we show that $p$ is an $r$-covering map. It is clear that $p$ is a graph homomorphism and $p|_{N(\alpha)}: N(\alpha) \rightarrow N(p(\alpha))$ is surjective for all $\alpha \in V(\tilde{G})$. By Lemma \ref{Lemma 6.1}, it suffices to show that $p|_{N_r(\alpha)}$ is injective for all $\alpha$. Let $\alpha \in V(\tilde{G})$, $\beta, \beta' \in N_r(\alpha)$, and suppose $p(\beta) = p(\beta')$. There are sequences $\gamma_0, \cdots, \gamma_r$, and $\gamma'_0,\cdots, \gamma'_r$ of vertices of $\tilde{G}$ such that $\gamma_0 = \gamma'_0 = \alpha$, $\gamma_r = \beta$, $\gamma'_r = \beta'$, and $(\gamma_{i-1},\gamma_i)$ and $(\gamma'_{i-1},\gamma'_i)$ belong to $E(\tilde{G})$ for $i=1,\cdots, r$. Let $\varphi$ be a representative of $\alpha$ and put $n = l(\varphi)$. Define the paths $\psi ,\psi': L_{n+r} \rightarrow G$ by $\psi |_{L_n} = \psi'|_{L_n} = \varphi$, $\psi(n+i)  = p(\gamma_i)$, and $\psi'(n+i) = p(\gamma'_i)$. Then we have that $\psi \in \beta$ and $\psi' \in \beta'$ (see (3) mentioned above). Since $\psi \simeq_r \psi'$, we have $\beta= \beta'$. Hence the map $p|_{N_r(\alpha)}$ is injective. Thus $p$ is an $r$-covering map.

It remains to be shown that $\tilde{G}$ is connected and $\pi_1^r(\tilde{G},\tilde{v})$ is trivial. Let $\alpha \in V(\tilde{G})$ and let $\varphi:L_n \rightarrow G$ be a representative of $\alpha$. Define the path $\tilde{\varphi}$ of $\tilde{G}$ by $\tilde{\varphi} (i) = [\varphi |_{L_i}]_r$. Note that the path $\tilde{\varphi}$ connects $\tilde{v}$ with $[\varphi]_r = \alpha$. Hence the graph $\tilde{G}$ is connected. Moreover, if $\alpha \in p_* (\pi_1^r(\tilde{G},\tilde{v})) \subset \pi_1^r(G,v)$, then we have $[*_v]_r = \tilde{v} = \tilde{\varphi}(n) = \alpha$. This implies that $p_* (\pi_1^r(\tilde{G},\tilde{v}))$ is trivial. It follows from Corollary \ref{Corollary 6.4} that $\pi_1^r(\tilde{G},\tilde{v})$ is trivial.
\end{proof}

Next we consider the relationship between $r$-covering maps and group actions.

Throughout the section, all group actions are assumed to be from the right. Let $\Gamma$ be a group and consider a $\Gamma$-action on the graph $G$. Namely, $\Gamma$ acts on the vertex set $V(G)$ from the right, and for every $\gamma \in \Gamma$, the map $V(G) \rightarrow V(G)$, $v \mapsto v\gamma$ is a graph homomorphism. Define the graph $G/\Gamma$ as follows. The vertex set is the orbit set $V(G) / \Gamma$, and two orbits $\alpha$ and $\beta$ are adjacent if and only if $(\alpha \times \beta) \cap E(G) \neq \emptyset$. Clearly, the quotient map $V(G) \rightarrow V(G/\Gamma)$ is a graph homomorphism from $G$ to $G/\Gamma$. The $\Gamma$-action is an {\it $r$-covering action} if $1 \neq \gamma \in \Gamma$ implies $N_r(v) \cap N_r(v \gamma) = \emptyset$ for every vertex $v$ of $G$.

\begin{prop} \label{Proposition 6.7}
Let $G$ be a graph having no isolated vertices, and let $\Gamma$ be a group. Suppose that $\Gamma$ acts on $G$, and the action is free as a set action on $V(G)$. Then the graph homomorphism $p:G \rightarrow G/\Gamma$ is an $r$-covering if and only if the $\Gamma$-action is an $r$-covering action.
\end{prop}
\begin{proof}
Suppose that the $\Gamma$-action is an $r$-covering action. By Lemma \ref{Lemma 6.1}, it suffices to show that $p|_{N(v)} : N(v) \rightarrow N(p(v))$ is surjective and $p|_{N_r(v)} : N_r(v) \rightarrow N_r(p(v))$ is injective for every $v \in V(G)$.

We show that $p|_{N(v)}$ is surjective. Let $\alpha$ be an element of $N(p(v))$. By the definition of $G/\Gamma$, there are an element $\gamma$ of $\Gamma$ and a representative $w$ of $\alpha$ such that $(v \gamma, w) \in E(G)$. Therefore we have $w \gamma^{-1} \in N(v)$ and $p( w\gamma^{-1}) = \alpha$. Hence $p|_{N(v)}$ is surjective.

We show that $p|_{N_r(v)}$ is injective. Let $w_0$ and $w_1$ be elements of $N_r(v)$ and suppose $p(w_0) = p(w_1)$. By the definition of $G/\Gamma$, there is $\gamma \in \Gamma$ with $w_0 \gamma = w_1$. Note
$$v \in N_r(w_0) \cap N_r(w_1) = N_r(w_0) \cap N_r(w_0 \gamma)$$
and hence we have $N_r(w_0) \cap N_r(w_0 \gamma) \neq \emptyset$. Since the $\Gamma$-action is an $r$-covering action, we have $\gamma = 1$ and $w_1 = w_0 1 = w_0$. Thus $p|_{N_r(v)} : N_r(v) \rightarrow N_r(p(v))$ is injective. Hence we have shown that the quotient map $p : G \rightarrow G / \Gamma$ is an $r$-covering map.

Next suppose that the $\Gamma$-action is free and the quotient map $p: G \rightarrow G/ \Gamma$ is an $r$-covering map. Let $v \in V(G)$, let $\gamma \in \Gamma$, and suppose $N_r(v) \cap N_r(v \gamma) \neq \emptyset$. Note that for an element $w \in N_r(v) \cap N_r(v \gamma)$, we have $v, v \gamma \in N_r(w)$ and $p(v) = p(v \gamma)$. Since $p$ is an $r$-covering map, we have $v = v \gamma$. Since the action is free, we have $\gamma = 1$. Thus  the action is an $r$-covering action.
\end{proof}

Let $p:(\tilde{G},\tilde{v}) \rightarrow (G,v)$ be the universal $r$-covering. We construct the bijection
$$\begin{CD}
\Phi : \pi_1^r(G,v) @>{\cong}>> p^{-1}(v)
\end{CD}$$
as follows. Let $\alpha$ be an element of $\pi_1^r(G,v)$ and let $\varphi$ be an element of $\alpha$. By (1) of Proposition \ref{Proposition 6.3}, there is a unique lift $\tilde{\varphi}$ starting from $\tilde{v}$. Let $\Phi([\varphi]_r)$ be the terminal point of $\tilde{\varphi}$. The independence of the choice of a representative $\varphi$ of $\alpha$ is deduced from (2) of Proposition \ref{Proposition 6.3}. The verification of the fact that $\Phi$ is bijective is straightforward and is left to the reader.

The above construction gives an alternative proof of $\pi_1^r(C_{2n+1}) \cong \Z$ with $r \leq n+1$ (Example \ref{Example 3.4}) without using Theorem \ref{Theorem 3.1}. Let $L$ be the graph defined by
$$V(L) = \Z,$$
$$E(L) = \{ (x,y) \; | \; |x-y| \leq 1\}.$$
Consider $0$ as a basepoint of $L$. Then it is easy to see that $\pi_1^r(L)$ is trivial. In fact, if $\varphi : L_n \rightarrow L$ is a loop of $L$ and $a \in V(L_n)$ is a point such that $|\varphi(a)|$ is maximum, then we have $\varphi(a-1) = \varphi(a+1)$. Moreover, the natural projection $L \rightarrow C_m$ is an $r$-covering map. Thus the above $\Phi$ gives a bijection
$$\Z \rightarrow \pi_1^r(C_{2n+1}), k \mapsto [\varphi]_r^k.$$
Here $\varphi : L_{2n+1} \rightarrow C_{2n+1}$ is the map $i \mapsto (i \; {\rm mod}. 2n+1)$. Note that $[\varphi]_r^2 = 1$ in $\pi_1^r(C_{2n+1})$ when $r > n+1$. Thus we also have $\pi_1^r(C_{2n+1}) = \Z / 2\Z$ when $r > n+1$. The computation of $\pi_1^r(C_{2n})$ is similarly obtained.

The following is the main result in this section.

\begin{thm} \label{Theorem 6.8}
Let $(G,v)$ be a based graph, and let $\Gamma$ be a subgroup of $\pi_1^r(G,v)$. Then there is a connected basepoint preserving $r$-covering $p_\Gamma:(G_\Gamma,v_\Gamma) \rightarrow (G,v)$ such that $p_{\Gamma*}(\pi_1^r(G_\Gamma,v_\Gamma)) = \Gamma$. Moreover, such an $r$-covering is unique up to isomorphisms.
\end{thm}
\begin{proof}
The uniqueness of $(G_\Gamma,v_\Gamma)$ follows from Lemma \ref{Lemma 6.5}. If $v$ is isolated, the proof is trivial. So we assume that $v$ is not isolated.

Consider the universal $r$-covering $(\tilde{G}, \tilde{v})$ constructed in the proof of Proposition \ref{Proposition 6.6}. Since $v$ is not isolated, $\tilde{G}$ has no isolated vertices. Define the $\pi_1^r(G,v)$-action on $\tilde{G}$ by
$$V(\tilde{G}) \times \pi_1^r(G,v) \rightarrow V(\tilde{G}), \; (\beta, \alpha) \mapsto \beta \cdot \alpha.$$
It is easy to see that for each element $\alpha$ of $\pi_1^r(G,v)$, the map $V(\tilde{G}) \rightarrow V(\tilde{G})$, $\beta \mapsto \beta \cdot \alpha$ is a graph homomorphism.

Let $\beta \in V(\tilde{G})$ and set $w = p(\beta)$. Then the orbit of the action is $\pi_1^r(G,v,w) = p^{-1}(w)$. Thus the induced map $V(\tilde{G} / \pi_1^r(G,v)) \rightarrow V(G)$ is bijective. It is easy to see that this is an isomorphism of graphs. Note that for $\beta \in V(\tilde{G})$ and $\alpha \in \pi_1^r(G,v)$, $\beta \cdot \alpha = \beta$ implies $\alpha = 1$ and hence this $\pi_1^r(G,v)$-action is free. Since the projection $p : \tilde{G} \rightarrow \tilde{G}/\Gamma \cong G$ is an $r$-covering map, we have that this $\pi_1^r(G,v)$-action is an $r$-covering action (see Proposition \ref{Proposition 6.7}).

Let $\Gamma$ be a subgroup of $\pi_1^r(G,v)$. Let $G_\Gamma = G / \Gamma$, $q: \tilde{G} \rightarrow G_\Gamma$ the quotient homomorphism, $p_\Gamma : G_\Gamma \rightarrow G$ the homomorphism such that $p_\Gamma \circ q = p$, and let $v_\Gamma = q(\tilde{v})$. By Proposition \ref{Proposition 6.7}, we have that $q$ is an $r$-covering map. By Lemma \ref{Lemma 6.2}, we have that $p_\Gamma$ is an $r$-covering map. To deduce $p_* \pi_1^r(G_\Gamma, v_\Gamma) = \Gamma$, it suffices to note the following commutative diagram:
$$\begin{CD}
\pi_1^r(G_\Gamma,v_\Gamma) @>{\cong}>> q^{-1}(v_\Gamma) @= \Gamma\\
@V{p_{\Gamma *}}VV @VVV @VVV\\
\pi_1^r(G,v) @>{\cong}>> p^{-1}(v) @= \pi_1^r(G,v),
\end{CD}$$
where the central and the right vertical arrows are inclusions. For the definitions of the left horizontal arrows, see the paragraph after Proposition \ref{Proposition 6.7}.
\end{proof}

Next we study the connected $r$-covering associated to the even part $\pi_1^r(G,v)_{ev}$ of $\pi_1^r(G,v)$. Let $G$ and $H$ be graphs. Define the {\it categorical product $G \times H$} by
$$V(G \times H) = V(G) \times V(H)$$
and
$$E(G \times H) = \{ ((x,y),(x',y')) \; | \; (x,x') \in E(G), \; (y,y') \in E(H)\}.$$
The {\it Kronecker double covering over $G$} (see \cite{IP}) is the 2nd projection $p : K_2 \times G \rightarrow G$. For a connected graph $G$, the Kronecker double covering over $G$ is connected if and only if $G$ is not bipartite. It is easy to show that the Kronecker double covering is an $r$-covering for every positive integer $r$.

\begin{prop} \label{Proposition Kronecker}
Let $(G,v)$ be a connected based graph and suppose that $G$ is not bipartite. Then the double covering associated to the even part of $\pi_1^r(G,v)$ is the Kronecker double covering $p : (K_2 \times G, (0,v)) \rightarrow (G,v)$ over $G$.
\end{prop}
\begin{proof}
Let $\varphi : L_n \rightarrow G$ be a loop of $(G,v)$. Then the lift of $\varphi$ with respect to $p : (K_2 \times G, (0,v)) \rightarrow (G,v)$ is the map $\tilde{\varphi} : L_n \rightarrow K_2 \times G$, $i \mapsto (i \; {\rm mod.} 2, \varphi(i))$. Note that $\tilde{\varphi}$ is a loop if and only if the length $n$ of $\varphi$ is even. Thus this propoerition follows from Corollary \ref{Corollary 6.4.1}.
\end{proof}

\begin{cor} \label{Corollary 6.21}
Let $(G,v)$ be a based graph. Then there is a 1-1 correspondence between connected based $(2r)$-coverings over $G$ and connected based $r$-coverings over $(K_2 \times G, (0,v))$.
\end{cor}
\begin{proof}
Note that when $G$ is bipartite, then $K_2 \times G$ is two copies of $G$. Thus it suffices to refer to Theorem \ref{Theorem 6.8}, Proposition \ref{Proposition Kronecker}, and Theorem \ref{Theorem 1.1}.
\end{proof}

The following proposition makes the correspondence of Corollary \ref{Corollary 6.21} apparent.

\begin{prop} \label{Proposition 6.20}
Let $p : G \rightarrow H$ be a $(2r)$-covering map. Then the map $p_* : N_r(G) \rightarrow N_r(H)$ induced by $p$ is a covering map.
\end{prop}
\begin{proof}
Let $w \in V(N_r(H))$. It suffices to show the following assertions:
\begin{itemize}
\item[(1)] If $v_1, v_2 \in p^{-1}(w)$ with $v_1 \neq v_2$, then ${\rm st}_{N_r(G)}(v_1) \cap {\rm st}_{N_r(G)}(v_2) = \emptyset$
\item[(2)]$p_*^{-1}({\rm st}_{N_r(H)}(w)) = \coprod_{v \in p^{-1}(w)} {\rm st}_{N_r(G)}(v)$
\item[(3)] For each $v \in p^{-1}(w)$, the simplicial map $p_* |_{{\rm st}_{N_r(G)}(v)} : {\rm st}_{N_r(G)}(v)\rightarrow {\rm st}_{N_r(H)}(w)$ is an isomorphism.
\end{itemize}

We first show (1). Suppose that there exists an element $v'$ of $V({\rm st}(v_1)) \cap V({\rm st}(v_2))$. Then $v_1, v_2 \in N_{2r}(v')$ and $p(v_1) = w = p(v_2)$. Since $p$ is a $(2r)$-covering, we have $v_1 = v_2$.

Next we show (2). Let $\sigma$ be a non-empty simplex of $p^{-1}_*({\rm st}_{N_r(H)}(w))$. Then we have $\sigma \in N_r(G)$ and $p(\sigma) \in {\rm st}_{N_r(H)}(w)$. Thus there exists a vertex $w'$ of $H$ such that $p(\sigma) \cup \{ w\} \subset N_r(w')$. For each element $x \in \sigma$, there is $v'_x$ such that $p(v'_x) = w'$. Then we have $v'_x = v'_y$ for every pair of elements $x$ and $y$ of $\sigma$. In fact, $p(v'_x) = w' = p(v'_y)$ and $v'_x, v'_y \in N_{2r}(x)$, where $x$ is a vertex of $G$ with $\sigma \subset N_r(x)$ . Set $v' = v'_x$. Then there is an element $v \in N_r(v')$ with $p(v) = w$. Since $\sigma \cup \{ v\} \subset N_r(v')$, we have that $\sigma \in {\rm st}(v)$. This completes the proof of the inclusion ``$\subset$'' in (2). The other direction is obvious.

Finally, we show (3). It is clear that $V({\rm st}(v)) = N_{2r}(v)$ and $V({\rm st}(w)) = N_r(w)$. Therefore $p|_{N_{2r}(v)} = p|_{{\rm st}(v)} : V({\rm st}(v)) \rightarrow V({\rm st}(w))$ is a bijection. Thus it suffices to show that the inverse $p|_{N_{2r}(v)}^{-1} : N_{2r}(w) \rightarrow N_{2r}(v)$ is a simplicial map.

Let $\sigma$ be a simplex of ${\rm st}(w)$ and set $\sigma' = (p|_{N_{2r}(v)})^{-1}(\sigma)$. There exists a vertex $w'$ of $H$ such that $\sigma \cup \{ w\} \subset N_r(w')$. Then there exists $v' \in N_r(v)$ with $p(v') = w'$. Let $\sigma$ be a subset of $N_r(v')$ such that $p(\sigma'') = \sigma$. Note that $\sigma'' \in {\rm st}(v)$ since $\sigma'' \cup \{ v\} \subset N_r(v')$. Since $\sigma', \sigma'' \subset N_{2r}(v)$ and $p(\sigma') = p(\sigma'')$, we have that $\sigma' = \sigma''$. This completes the proof.
\end{proof}

Now we consider the correspondence of Corollary \ref{Corollary 6.21}. Let $(G,v)$ be a based graph such that $v$ is not isolated. We write $N_r(G)_0$ to indicate the connected component of $N_r(G)$ containing $v$.

Suppose that $G$ is non-bipartite. It is clear that $N_r(K_2 \times G)$ is two copies of $N_r(G)$. Let $p : H \rightarrow K_2 \times G$ be a connected basepoint preserving $(2r)$-covering over $K_2 \times G$. Then $N_r(H)_0 \rightarrow N_r(K_2 \times G)_0 = N_r(G)_0$ is a covering space (Proposition \ref{Proposition 6.20}). Similarly, if $G$ is bipartite, then a connected basepoint preserving $(2r)$-covering $p : H \rightarrow G$, then $p_* : N_r(H)\rightarrow N_r(G)_0$ is a covering space. These are the correspondences of Corollary \ref{Corollary 6.21}.

We conclude this paper with a few remarks. Let $G$ and $H$ be connected non-bipartite graphs and suppose that the Kronecker double coverings over $G$ and $H$ are isomorphic. Then we have that $\pi_1^2(G)_{ev} \cong \pi_1^2(H)_{ev}$. By Theorem \ref{Theorem 1.1}, we have $\pi_1(N(G)) \cong \pi_1(N(H))$. In fact, the author showed that $K_2 \times G \cong K_2 \times H$ implies $N(G) \cong N(H)$ (see \cite{Matsushita}). 

It follows from Theorem \ref{Theorem 1.1}, Theorem \ref{Theorem 5.2}, and Theorem \ref{Theorem 6.8} if the neighborhood complex of $G$ is simply connected, then a connected 2-covering over $G$ is isomorphic to either $G$ or $K_2 \times G$. Examples of such graphs are given by the Kneser graphs $KG_{n,k}$ for $n > 2k +1$ (see Theorem \ref{Theorem 5.2}) or some of the stable Kneser graphs discussed in \cite{BL}. This phenomenon is quite different from the usual covering maps over graphs. In fact if a connected graph $G$ has an embedded cycle, then there are infinitely many coverings over $G$.

\vspace{2mm} \noindent {\bf Acknowledgements.} I wish to express my gratitude to Professor Toshitake Kohno for his indispensable advice and support. I thank Shouta Tounai for stimulated conversations. I thank the anonymous referees for helpful comments. I am supported by the Grant-in-Aid for Scientific Research (KAKENHI No. 25-4699 and 28-6304), the Grant-in-Aid for JSPS fellows, and the Program for Leading Graduate Schools, MEXT, Japan.


\begin{thebibliography}{99}

\bibitem{BL} A. Bj\"orner, M. D. Longueville, {\it Neighborhood complexes of stable Kneser graphs}, {\bf 23} (1) 23-34 (2003).

\bibitem{IP} W. Imrich, T. Pisanski, {\it Multiple Kronecker covering graphs}, European J. Combin. {\bf 29} 1116-1122 (2008).

\bibitem{Kozlov} D. N. Kozlov. {\it Combinatorial algebraic topology.} Algorithms and Computation in Mathematics. Vol. 21 Springer, Berlin (2008).

\bibitem{Lovasz} L. Lov\'asz, {\it Kneser's conjecture, chromatic number, and homotopy}, J. Combin. Ser. A {\bf 25} (3) 319-324 (1978).

\bibitem{MZ} J. Matou$\check{\rm s}$ek, G. M. Ziegler, {\it Topological lower bounds for the chromatic number: A hierarchy}, Jahresbericht der Deutschen Mathematiker-Vereinigung {\bf 106} (2004) no. 2, 71-90 

\bibitem{Spanier} E. H. Spanier, {\it Algebraic topology}, Springer-Verlag (1966).

\bibitem{Matsushita} T. Matsushita, {\it Box complexes and Kronecker double coverings of graphs}, arXiv:1404.1549

\end{thebibliography}
\end{document}